\documentclass[oneside]{amsart}
\usepackage{float}
\usepackage{amsmath}
\usepackage{graphicx}
\usepackage{amsfonts}
\usepackage{amssymb}
\usepackage[cmtip,all]{xy}

\theoremstyle{plain}
\newtheorem{Theorem}{Theorem}[section]

\newtheorem{Lemma}[Theorem]{Lemma}
\newtheorem{Proposition}[Theorem]{Proposition}
\newtheorem{Corollary}[Theorem]{Corollary}

\theoremstyle{remark}

\newtheorem{Remark}[Theorem]{Remark}

\theoremstyle{definition}

\numberwithin{equation}{section}

\begin{document}
\title[Multiplicity-Free $K_{\mathbb{C}}$-Orbits]
{On a Class of Multiplicity-Free Nilpotent $K_{\mathbb{C}}$-Orbits}
\author{B. Binegar}
\address{Department of Mathematics\\
Oklahoma State University\\
Stillwater, Oklahoma 74078}
\email{binegar@okstate.edu}
\date{October 4, 2007}
\thanks{The author gratefully acknowledges discussions with Don King
and Kyo Nishiyama at the Snowbird Conference on Representations of Real Reductive
Groups, June, 2006. He also thanks the National Science Foundation through the Atlas 
for Lie Groups FRG (DMS 0554278) for support.}
\maketitle
\begin{abstract} Let $G$ be a real, connected, noncompact, semisimple Lie group,
let $K_{\mathbb{C}}$ be the complexification of a maximal compact subgroup $K$ of $G$, and let $\frak{g}=\frak{k}+\frak{p}$ be
the corresponding Cartan decomposition of the complexified Lie algebra of $G$.
Sequences of strongly orthogonal noncompact weights are constructed and
classified for each real noncompact simple Lie group of classical type. We
show that for each partial subsequence $\left\{\gamma_{1},\ldots, \gamma
_{i}\right\}$ there is a corresponding family of nilpotent $K_{\mathbb{C}}$-orbits
in $\frak{p}$, ordered by
inclusion and such that the representation of $K$ on the ring of regular
functions on each orbit is multiplicity-free. The $K$-types of regular functions
on the orbits and the regular functions on their closures are both
explicitly identified and demonstrated to coincide, 
with one exception in the Hermitian symmetric case. 
The classification presented also includes
the specification of a base point for each orbit and exhibits a corresponding
system of restricted roots with multiplicities. A formula for the
leading term of the Hilbert polynomials corresponding to these orbits is given.
This formula, together with the restricted root data, allows the determination
of the dimensions of these orbits and the algebraic-geometric degree of their
closures. In an appendix, the location of these orbits within D. King's
classification of spherical nilpotent orbits in complex symmetric spaces
is depicted via signed partitions and Hasse diagrams.
\end{abstract}

\section{Introduction}

An action of an algebraic reductive group $G$ on an affine variety $M$ is
called \emph{multiplicity-free} if the multiplicity of any particular
irreducible representation of $G$ in the space $\mathbb{C}\left[  M\right]  $
of regular functions on $M$ is at most one. In \cite{kac}, Kac provides a complete
list of multiplicity-free actions for the situation where $G$ is a connected
reductive algebraic group and $M$ is a finite-dimensional vector space upon
which $G$ acts by an irreducible representation. We remark that Kac initiated
this classification in order to understand the possibilities for the $G_{0}$-orbits 
in $\frak{g}_{i}$, where $\frak{g}_{i}$ is $i^{th}$ homogeneous component of 
a $\mathbb{Z}$-graded semisimple Lie algebra and $G_{0}$ is the adjoint group 
of $\frak{g}_{0}$. 

In \cite{kato-ochiai}, Kato and Ochiai develop a formula for the algebraic-geometric degree
of a multiplicity-free $G$-variety $Y$ in the situation where $G$ is a
connected reductive complex algebraic group, and $Y$ is a closed $G$-stable subset of a 
finite-dimensional vector space $V$ carrying an irreducible, multiplicity-free representation of $G$
and such that the image of $G$ in $GL\left(  V\right)$
contains all nonzero scalar matrices. Kato and
Ochiai then proceed to explicitly evaluate their formula for the case when $V$ is
the holomorphic tangent space of a Hermitian symmetric space $G/K$ regarded as 
a representation of the complexification $K_{\mathbb{C}}$ of $K$. In this
last situation, there exists a set of
linearly independent dominant weights $\left\{  \varphi_{1},\ldots,\varphi
_{i}\right\}  $ so that
\[
\mathbb{C}\left[  Y\right]  \cong\bigoplus_{m\in\mathbb{N}^{i}}V_{m_{1}%
\varphi_{1}  +\cdots +m_{i}\varphi_{i}  }%
\]
where $V_{m_{1}\varphi_{1}+\cdots +m_{i}\varphi_{i} }$ denotes the irreducible representation of $K$ of highest 
weight $m_{1}\varphi_{1}  +\cdots +m_{i}\varphi_{i}$ and the sum is over all $m$-tuples of non-negative integers 
$(m_{1},\ldots,m_{i})$.
Moreover, in the Hermitian symmetric situation, there is a natural way of constructing the
weights $\varphi_{j}$, $j=1,\ldots,i$ from a subsequence $\left\{\gamma_{1},\ldots,\gamma_{i}\right\}$
of a Harish-Chandra sequence 
$\left\{  \gamma_{1},\ldots,\gamma_{n}\right\}$ of strongly orthogonal non-compact roots, as well as an explicit
accounting of the roots that contribute, via the Weyl dimension formula, to
the degree of the orbit. It happens that the contributing (restrictions of) 
positive roots break up into two disjoint subsets%
\begin{align*}
\Delta_{short}^{+}  &  =\left\{  \frac{1}{2}\gamma_{j}\mid 1\leq j\leq
i\right\}  \text{ with a common multiplicity }r\\
\Delta_{long}^{+}  &  =\left\{  \frac{1}{2}\gamma_{j}-\frac{1}{2}\gamma
_{k}\mid 1 \leq j<k \leq i\right\} \text{ with a common multiplicity } k  
\end{align*}
These circumstances allow Kato and Ochiai to reduce the problem of determining the
algebraic-geometric degree of $Y$ to an application
of the Selberg integral formula (\cite{selberg}).

From the ``orbit philosophy'' point of view in representation theory, there are two other 
especially important, general cases 
of multiplicity-free actions: the case when $M$ is nilpotent $Ad(\frak{g})$-orbit in the Lie algebra 
of a complex semisimple Lie algebra $\frak{g}$ for which a Borel subgroup of 
$Ad(\frak{g})$ has a dense orbit, and the case when $M$ is an
an irreducible component of the associated
variety of a multiplicity-free $\left(  \frak{g},K_{\mathbb{C}}\right)  $-module. The orbits in 
the first case are called \emph{spherical} nilpotent orbits and these have been studied
and classified by Panyushev \cite{panyushev}. See also \cite{KY}, where spherical nilpotent orbits for
a complex Lie algebra are realized within the secant variety attached to the adjoint variety of a simple
complex Lie algebra.  

The associated varieties in the 
second case correspond to multiplicity-free $K_{\mathbb{C}}$-orbits in $\mathcal{N}_{\frak{p}}$,
the nilpotent cone in $\left(\frak{g} \backslash \frak{k} \right)^{\ast} \cong \frak{p}$. Such
orbits are referred to as spherical nilpotent orbits for the symmetric pair
$\left(\frak{g},\frak{k}\right)$. These have been classified by D. King (\cite{king}). We remark
that the varieties $Y$ studied by Kato and Ochiai can be viewed as 
a special cases (the Hermitian symmetric cases) of a spherical nilpotent orbit for
a symmetric pair. We note further the papers \cite{N1}, \cite{N2}, \cite{NO}, \cite{NOT}, \cite{noz};
wherein the associated varieties of singular unitary representations attached 
to certain dual pairs are shown to be multiplicity free. In the last three papers, integral 
formulas for the Bernstein degrees of the representations are also developed and in some cases
explicitly computed. In particular, in \cite{NOT} it is observed that the explicit formulas
for Bernstein degrees so obtained coincide with the classical Giambelli formulas for
the degrees of determinantal varieties. In fact, such integral formulas for degrees are common to 
spherical varieties in general ( \cite{Br1}, \cite{Br2} ). 

In this paper, we reverse-engineer the results of Kato and Ochiai to obtain
a construction of a family of multiplicity-free $K_{\mathbb{C}}$-orbits 
in $\mathcal{N}_{\frak{p}}$ that is applicable for any noncompact semisimple Lie algebra $\frak{g}$. 
However, instead of starting with $K_{\mathbb{C}}$-orbits known
to be multiplicity free, and looking for an
associated sequence of strongly orthogonal noncompact roots; we proceed as follows:

\begin{enumerate}
\item  In the context of an arbitrary connected
noncompact real semisimple Lie group $G$ we introduce an algorithm for constructing
sequences $\left\{  \gamma_{1},\ldots,\gamma_{n}\right\}  $ of strongly
orthogonal noncompact weights.

\item  We then attach to each subsequence $\left\{\gamma_{1},\ldots,\gamma_{i}\right\}$
a certain nilpotent element $Y_{i}$ of $\frak{p}$, and set 
$\mathcal{O}_{i} = K_{\mathbb{C}}\cdot Y_{i}$.
We show that the closure $\overline{\mathcal{O}_{i}}$ of each such orbit is
multiplicity-free, and we explicitly identify the $K$-types of regular functions on
the closure $\overline{\mathcal{O}_{i}}$ of $\mathcal{O}_{i}$ as
\[
\mathbb{C}\left[  \overline{\mathcal{O}_{i}}\right]  \cong\bigoplus
V_{a_{1}\gamma_{1}+\cdots+a_{i}\gamma_{i}}%
\]
where the sum is over the $a_{i}\in\mathbb{N}$ such that $a_{1}\geq
a_{2}\geq\cdots\geq a_{i}\geq 0$. (However, when the restricted root system is type $D_n$, and 
$\left\{ \gamma_{1} , \ldots, \gamma_{n} \right\}$ is a sequence of maximal length, the bound
on the last coefficient is actually $\left| a_{n} \right| \geq 0.$) 

\item We observe that the degree of homogeneity of 
a polynomial in
$V_{a_{1}\gamma_{1}+\cdots+a_{i}\gamma_{i}}$ is $\sum_{j=1}^{i}a_{j}$, and thereby
reproduce the canonical filtration of $\mathbb{C}\left[  \overline{\mathcal{O}_{i}}\right]  $
by degree by setting
\begin{equation}
\mathbb{C}\left[  \overline{\mathcal{O}_{i}}\right]  _{\ell}\cong
\bigoplus_{\substack{a_{1}\geq a_{2}\geq\cdots\geq a_{i}\geq0\\\sum a_{j}%
\leq\ell}}V_{a_{1}\gamma_{1}+\cdots+a_{i}\gamma_{i}} \label{fil}%
\end{equation}
Using the restricted root data obtained in \S 2 and the Weyl dimension formula, 
we are then able to calculate the leading term of the corresponding Hilbert polynomial 
and thereby obtain formulas for the dimension and algebraic-geometric degree of
$\overline{\mathcal{O}_i}$ in the classical cases. 
\end{enumerate}

The organization of this paper is as follows. In \S 2 we define certain sequences
of strongly orthogonal noncompact weights. These sequences will provide the basic substratum
upon which everything else is pinned.  Table 1 in
that section lists, for each real classical Lie group, the sequences of strongly orthogonal
noncompact weights of maximal length and the form of their restricted root systems (as defined
in that section). 

In \S 3 we attach to each sequence of strongly orthogonal noncompact weights $\Gamma = \left\{\gamma_{1},\ldots,\gamma_{n}\right\}$
a corresponding sequence $\left\{x_{i},h_{i},y_{i}\right\}$, $i=1,\ldots,n$ of mutually centralizing normal $S$-triples. These 
in turn allow us to construct the ``telescoping'' sequences of $K_{\mathbb{C}}$-orbits
$\mathcal{O}_{1} \subset \mathcal{O}_{2} \subset \cdots \subset \mathcal{O}_{n}$ which will be the principal objects of study for the
rest of the paper.  We show that each orbit $\mathcal{O}_{i}$ is
multiplicity-free and determine the $K$-type decompositions of the rings of regular functions on $\mathcal{O}_{i}$ and its closure.

We conclude \S 3 with two remarks; the first indicating where our family of nilpotent $K_{\mathbb{C}}$-orbits sits within D. King's
\cite{king} classification of nilpotent orbits for classical symmetric pairs. The second remark sketches our plan to attach to
such a family of $K_{\mathbb{C}}$-orbits a corresponding family of unipotent representations.  We note that effectively this has already been 
achieved by Sahi in the situation where $G$ is the conformal group of a Euclidean \cite{sahi0} or non-Euclidean \cite{sahi}
real simple Jordan algebra.  We show in \S3.1.2 how one can recover, in the context of
an arbitrary connected semisimple Lie group, nearly all of the structural niceties employed by 
Sahi in \cite{sahi} to bring to light families of unitarizable unipotent representations residing within  
families of degenerate principal series representations attached to corresponding families of nilpotent
$K_{\mathbb{C}}$-orbits.  

In \S 4 we utilize the $K$-type decompositions determined in \S 3 and along with the forms of the restricted root systems
given in Table 1 to obtain closed formulas for the dimension and algebraic-geometric degrees of the closures of the orbits. 
We thereby produce analogs of the formulas of Kato and Ochiai in the general setting of noncompact classical Lie groups.  

The multiplicity-free $K_{\mathbb{C}}$-orbits of real classical noncompact groups that we obtain in this paper
all lie within the King's classification \cite{king} of spherical nilpotent orbits for
symmetric pairs (and we hereby apologize for adopting a nomenclature that might suggest otherwise).
In an appendix, we illustrate via Hasse diagrams how our orbits are situated amongst
the other orbits in King's classification. We remark that the $K$-type decompositions for \emph{all} the 
spherical orbits of the symmetric pairs $\left( U(p,p)/U(p)\times U(p) \right)$ 
have recently been obtained
by K. Nishiyama (\cite{N2}) using dual pair methods, while for the same symmetric pairs, our method yields only the 
spherical orbits that reside along the outer edges of the corresponding Hasse diagram. 

\section{Sequences of Strongly Orthogonal Noncompact Weights}

Let $G$ be a connected noncompact real semisimple Lie group. Let $K$ be a
maximal compact subgroup, $\theta$ the corresponding Cartan involution and
$\frak{g}=\frak{k}+\frak{p}$, the corresponding Cartan decomposition of the
complexification of the Lie algebra of $G$. Choose a Cartan subalgebra
$\frak{t}$ of $\frak{k}$, and extend it to a $\theta$-stable Cartan subalgebra
$\frak{h}=\frak{t}+\frak{a}$ of $\frak{g}$. Choose a positive system
$\Delta^{+}\left(  \frak{t};\frak{k}\right)  $ for $\Delta\left(
\frak{t};\frak{k}\right)  $ and extend it to a positive system $\Delta
^{+}\left(  \frak{h};\frak{g}\right)  $ of $\Delta\left(  \frak{h}%
;\frak{g}\right)  $ in such a way that%
\[
\left.  \alpha\right|  _{\frak{t}}\in\Delta^{+}\left(  \frak{t};\frak{k}%
\right)  \quad\Longrightarrow\quad\alpha\in\Delta^{+}\left(  \frak{h}%
;\frak{g}\right) \quad .
\]
Let $\widetilde{\beta}$ be a highest weight of an irreducible representation of $K$ on
$\frak{p}$. We remark that $\beta$ is unique when $\frak{g}$ is simple not of Hermitian
type.  In the simple Hermitian symmetric case, where $\frak{p}$
decomposes into a sum of two irreducibles, $\frak{p}=\frak{p}_{+}+\frak{p}%
_{-}$, and one can take $\beta$ to be the highest weight of the representation of $K$
on  $\frak{p}_{+}$ or $\frak{p}_{-}$. We now construct sequences $\left\{  \gamma
_{1},\ldots,\gamma_{n}\right\}  $ of strongly orthogonal noncompact weights as follows.

\begin{itemize}
\item We set $\gamma_{1}=\widetilde{\beta}$;

\item $\gamma_{i+1}$ is determined from its predecessors $\left\{  \gamma_{1},\ldots,\gamma
_{i}\right\}  $ by the requirements

\begin{itemize}
\item [(i)]$\gamma_{i+1}$ is in the orbit of $\widetilde{\beta}$ under the
action of the Weyl group of $K$.

\item[(ii)] $\gamma_{i+1}$ is strongly orthogonal in $\frak{g}$ to each
$\gamma_{j}$ for $j=1,,\ldots,i$ (meaning there is no compact or noncompact
weight vector of weight $\gamma_{i+1} \pm \gamma_{j}$ for $j=1,\ldots,i$.)

\item[(iii)] $\omega_{i+1}=\sum_{j=1}^{i+1}\gamma_{j}\in\frak{t}^{\ast}$ is dominant.
\end{itemize}
\end{itemize}
Of course, since $\dim \frak{p}$ is finite, this constructive process will eventually 
terminate. It turns out that, almost always, the maximal length of such a sequence 
is equal to the lesser of the rank of $K$ and the real rank of $G$. (See the remarks following
 Table 1.)  
\bigskip

In Table 1 below we tabulate, for each real classical noncompact Lie group of real rank
$\geq 2$, 
sequences $\Gamma = \left\{\gamma_{1},\ldots,\gamma_{n} \right\}$ of maximal length. 
(The real rank one cases are excluded simply by virtue of their triviality: 
in these cases $\Gamma = \left\{\widetilde{\beta} \right\}$.) 

We also provide in the table
the form of the \emph{restricted root systems} for $\Gamma$. This restricted root system is
defined as follows. For each noncompact weight $\gamma_{i} \in \Gamma$ we can choose a representative
nilpotent element $x_i$ in $\frak{p}_{\gamma_i}$, and then via a standard construction, a
normal $S$-triple $\left\{ x_{i}, h_{i} , y_{i}\right\} $ where $y_{i}\in \frak{p}_{-\gamma_{i}}$,
$h_{i} \in \frak{t}$ and
\[
\left[ x_{i} ,y_{i} \right] = h_{i} \quad, \quad \left[h_{i},x_{i}\right] = 2x_{i} \quad , 
\quad \left[ h_{i} , y_{i} \right] = -2y_{i} \quad .
\]
Set $\frak{t}_{1} =
span_{\mathbb{C}} \left( h_{1}, \ldots , h_{n} \right) \subset \frak{k}$. The restricted root
system $\Sigma$ corresponding to $\Gamma$ is the set of $\frak{t}_{1}$-roots in $\frak{k}$. In the
table, the form of a restricted root system $\Sigma$ is indicated follows:%
\begin{equation}
\Sigma=\left(  a_{n}\right)  ^{m_{a}}\left(  A_{n}\right)  ^{m_{A}}\left(
b_{n}\right)  ^{m_{b}}\left(  C_{n}\right)  ^{m_{C}}\left(  d_{n}\right)
^{m_{d}} \left( a_{11,+} \right)^{m_{+}} \left( a_{11,-} \right)^{m_{-}} \label{0}%
\end{equation}
means that the set of positive roots in $\Sigma$ consists of roots of the form

\begin{itemize}
\item  $a_{n} = \left\{  \frac{1}{2}\gamma_{i}-\frac{1}{2}\gamma
_{j}\mid1\leq i<j\leq n\right\}  $, each occurring with multiplicity $m_{a}$;

\item  $A_{n} = \left\{  \gamma_{i}-\gamma_{j}\mid1\leq i<j\leq
n\right\}  $, each occurring with multiplicity $m_{A}$;

\item  $b_{n}=\left\{  \frac{1}{2}\gamma_{i}\mid1\leq i\leq
n\right\}  $, each occurring with multiplicity $m_{b}$;

\item  $C_{n} = \left\{  \gamma_{i}\mid1\leq i\leq n\right\}  $, each
occurring with multiplicity $m_{C}$; 

\item  $d_{n} =\left\{  \frac{1}{2}\gamma_{i}\pm\frac{1}{2}%
\gamma_{j}\mid1\leq i<j\leq n\right\}  $, each occurring with multiplicity
$m_{d}$ ;

\item $a_{11,+}=\left\{ \pm \left(\frac{1}{2}\gamma_{1} +  \frac{1}{2}\gamma_{2}
\right) \right\}$, each occuring with multiplicity $m_{+}$; and 

\item $a_{11,-}=\left\{ \pm \left(\frac{1}{2}\gamma_{1} -  \frac{1}{2}\gamma_{2}
\right) \right\}$, each occuring with multiplicity $m_{+}$. 
\end{itemize}

We remark that $m_{a}\ne 0$ or $m_{A}\neq 0$ only in the Hermitian symmetric
case, and in this case $m_{d}=0$.\footnote{Note also that our mnemonic 
notation for restricted root systems is a bit misleading for the types 
$a_{n}$ and $A_{n}$ since these are actually root systems of Cartan type 
$A$ and rank $n-1$.}
We specify in Table 1 the non-compact weights $\gamma _{i}$ in terms of
a basis of fundamental weights of the semisimple part $\left[ K,K\right] $
of $K$ and the conventions of Bourbaki (\cite{bourbaki}).
When $\left[ K,K\right] $ has two factors, say for rank $r$ and $s$, we
denote by $\omega _{1},\ldots ,\omega _{r}$ a basis (\`{a} la Bourbaki) for
fundamental weights for the first factor, and $\omega _{r+1},\ldots ,\omega
_{r+s}$ a basis of fundamental weights for the second factor.

When $\left[ K,K\right] $ has an $SO\left( n\right) $ factor several
idiosyncrasies occur which we shall now describe in detail. First of all, we
have to deal with the fact that $SO\left( n\right) \sim D_{%
\left[ \frac{n}{2}\right] }$ when $n$ is even and $SO\left( n\right) \sim B_{%
\left[ \frac{n}{2}\right] }$ when $n$ is odd.  It also turns out that, for
even $n$, we have two different ways of terminating maximal sequences of
strongly orthogonal noncompact weights (corresponding to the outer
automorphism of $D_{n}$). We shall employ the following shorthand to deal
efficiently these variations. Let $\sigma _{n,\pm }$ denote the following
sequences of weights (of $SO\left( n \right)$). 
{\small 
\begin{equation*}
\sigma _{n,\pm }=\left\{ 
\begin{array}{ll}
\omega _{1}+\omega _{2} & \text{if }n=4 \quad , \\ 
\omega _{1},\;\omega _{2}-\omega _{1},\;\ldots ,\;
\omega _{k-1} +\omega _{k}-\omega _{k-2},\;\pm
\omega _{k}\mp \omega_{k-1} & \text{if } n=2k>4  \quad , \\ 
\omega _{1},\;\omega_{2}-\omega_{1},\;\ldots ,\;-\omega _{k-2}
+\omega _{k-1},\;2\omega _{k}-\omega _{k-1} & \text{if }n=2k+1\quad.
\end{array}
\right. 
\end{equation*}
}
We indicate by $2\sigma _{n,\pm }$, the sequences $2\omega _{1},\;2\omega
_{2}-2\omega _{1},\;\ldots $ \ , etc which occur in the case of $%
SL\left( n,\mathbb{R}\right) $.

To describe the sequences for $SO\left( p,q\right) $,  $p\leq q$, we first
denote by $\sigma _{p,q,\pm ,\pm }$ the sequence of $\left[ \frac{p}{2}%
\right] $ weights for $SO\left( q\right) \times SO\left( q\right) $ obtained
by adding to each element of the sequence $\sigma _{p,\pm }$ the
corresponding element in the sequence $\sigma _{q,\pm }$. Secondly, we
denote by $\tau _{p,q,\pm }$ the two-element sequences
\begin{equation*}
\tau _{p,q,\pm }=\omega _{1}+\omega _{\left[ \frac{p}{2}\right] +1}\ ,\ \pm
\omega _{i}\mp \omega _{\left[ \frac{p}{2}\right] +1} \quad .
\end{equation*}
The sequences of noncompact weights for \ $SO\left( p,q\right) $ will then
consist of the sequences$\ $ $\sigma _{p,q,\pm ,\pm }$ and $\tau _{p,q,,\pm }
$. Depending on the parities of $p$ and $q$, in the $SO\left( p,q\right) $
case, $3<p\leq q$,  there can be as many as six different sequences of
noncompact weights, or as few as three.

\begin{table}\caption{}
\scriptsize 
\begin{equation*}
\begin{array}{|c|c|c|l|}
\hline
G & K & \Sigma  & \qquad \qquad \Gamma  \\ \hline
SL\left( n,\mathbb{R}\right)  & SO\left( n\right)  & 
\begin{array}{ll}
 \left( d_{\left[ n/2 \right] }\right) ^{1} & \text{if } n \text{ is even}\\
 \left( b_{\left[ n/2 \right] }\right) ^{1} \left( d_{\left[ n/2 \right] }\right) ^{1} & \text{if } n \text{ is odd}
\end{array}
 & 2\sigma _{[n/2],\pm } \\ \hline
SL\left( n,\mathbb{H}\right)  & Sp\left( n\right)  & \left( C_{\left[ n/2%
\right] }\right) ^{3}\left( d_{\left[ n/2\right] }\right) ^{4} & 
\begin{array}{l}
\gamma _{1}=\omega _{2} \\ 
\gamma _{i}=\omega _{2i}-\omega _{2i-2} \\ 
\gamma _{\left[ n/2\right] }=\left\{ 
\begin{array}{l}
-\omega _{n-2}+\omega _{n}\text{ if }n\text{ is even} \\ 
-\omega _{n-3}+\omega _{n-1}\text{ if }n\text{ is odd}
\end{array}
\right. 
\end{array}
\\ \hline
\begin{array}{l}
SU\left( p,q\right) \  \\ 
\ (2\le p\leq q)
\end{array}
& S\left( U\left( p\right) \times U\left( q\right) \right)  & 
\left( a_{p}\right) ^{2}\left(b_{p}\right) ^{q-p} & 
\begin{array}{l}
\gamma _{1}=\omega _{1}+\omega _{p+q-1} \\ 
\gamma _{i}=-\omega _{i-1}+\omega _{i}+\omega _{p+q-i-1}-\omega _{p+q-i} \\ 
\gamma _{p}=-\omega _{p-1}+\omega _{q}-\omega _{q-1}
\end{array}
\\ \hline
\begin{array}{l}
SO\left( 2,q\right)  \\ 
\left( q>2\right) 
\end{array}
& S\left( O\left( 2\right) \times O\left( q\right) \right)  & \left(
A_{2}\right) ^{q-2} & 
\begin{array}{l}
\gamma _{1}=\omega _{1} \\ 
\gamma _{2}=-\omega _{1}
\end{array}
\\ 
\hline
\begin{array}{l}
SO\left( p,q\right), I  \\ 
\left( 2<p\leq q\right) 
\end{array}
& S\left( O\left( p\right) \times O\left( q\right) \right)  &
\left(b_{
\left[ \frac{p}{2}\right] }\right) ^{q-p+2\delta_p}\left( d_{\left[ \frac{p}{2}\right]
}\right) ^{2} 
& 
\sigma _{p,q,\pm,\pm} \\
\hline
\begin{array}{l}
SO\left( p,q\right), II  \\ 
\left( 2<p\leq q\right) 
\end{array}
& S\left( O\left( p\right) \times O\left( q\right) \right)  &
\left(a_{11,\pm}\right)^{p-2}\left( a_{11,\mp} \right)^{q-2}
& 
\tau_{p,q,\pm} 
\\ 
\hline
SO^{\ast }\left( 2n\right)  & U\left( n\right)  & \left( b_{\left[ \frac{n}{2%
}\right] }\right) ^{2}\left( d_{\left[ \frac{n}{2}\right] }\right) ^{4} & 
\begin{array}{l}
\gamma _{1}=\omega _{1} \\ 
\gamma _{i}=-\omega _{2\left( i-1\right) }+\omega _{2i} \\ 
\gamma _{\left[ \frac{n}{2}\right] }=\left\{ 
\begin{array}{l}
-\omega _{n-2}+\omega _{n}\quad \text{ if }n\text{ is even} \\ 
-\omega _{n-3}+\omega _{n-1}\quad \text{if }n\text{ is odd}
\end{array}
\right. 
\end{array}
\\ \hline
Sp\left( n,\mathbb{R}\right)  & U\left( n\right)  & \left( a_{n}\right) ^{1}
& 
\begin{array}{l}
\gamma _{1}=2\omega _{1} \\ 
\gamma _{i}=-2\omega _{i-1}+2\omega _{i} \\ 
\gamma _{n}=-2\omega _{n-1}
\end{array}
\\ 
\hline
\begin{array}{l}
Sp\left( p,q\right)  \\ 
\left( p\leq q\right) 
\end{array}
& Sp\left( p\right) \times Sp\left( q\right)  
&\left( b_{p}\right) 
^{2\left( q-p\right) } \left(C_{p}\right)^{2} \left( d_{p}\right) ^{2}
& 
\begin{array}{l}
\gamma _{1}=\omega _{1}+\omega _{p+1} \\ 
\gamma _{i}=-\omega _{i-1}+\omega _{i}-\omega _{p+i}+\omega _{p+i+1} \\ 
\gamma _{p}=-\omega _{p-1}+\omega _{p}-\omega _{2p}+\omega _{2p+1} \\
\end{array}
\\ \hline
\end{array}
\end{equation*}

(The term $\delta_{p}$ that appears in the exponent of $b_{\lbrack p/2 \rbrack}$ for
type $SO(p,q)$ is integer remainder of $p$ when divided by $2$.)  

\normalsize
\end{table}

\subsection{Remarks}
 
\subsubsection{}  One could consider relaxing the requirement that each $\gamma_{i}$
lie in the $K$-Weyl orbit of the highest noncompact weight by instead
stipulating that each $\gamma_{i}$ is a weight of the representation
of $K$ on $\frak{p}$. This leads to more sequences of strongly orthogonal
noncompact weights, but it seems that the sequences don't get any longer and, moreover,
our method of identifying the $K$-types supported on the closures of the orbits
is not applicable for such sequences.
In the Hermitian symmetric case, where $\frak{p}$ is a direct sum of two
irreducible representations of $K$, one could consider utilizing weights from both summands
to form strongly orthogonal sequences of noncompact weights. This does lead to additional 
long sequences of strongly orthogonal noncompact weights and, as K. Nishiyama
has pointed out to us,  the corresponding sequences of $K_{\mathbb{C}}$-orbits may actually 
exhaust the spherical nilpotent orbits for Hermitian symmetric pairs.  However, for such sequences
it is also difficult to identify exactly which $K$-types appear in the
ring of regular functions on the closures of the corresponding $K_{\mathbb{C}}$-orbits.  

\subsubsection{}  In all but the case of $SU^{\ast}\left(  2n\right)  $ the maximal
number of elements in a sequence of strongly orthogonal noncompact weights is
equal to $\min\left(  rank\left(  G/K\right)  ,rank\left(  K\right)  \right)
$. This suggests a connection with the maximal number of commuting
$\frak{sl}\left(  2,\mathbb{R}\right)  $ subalgebras of $\frak{g}_{\mathbb{R}}=Lie_{\mathbb{R}}(G)$.

Indeed, to each $\gamma_{i}\in\left[  \gamma_{1},\ldots,\gamma_{n}\right]  $
we have an associated normal triple $\left\{  x_{i},h_{i},y_{i}\right\}  $.
In fact, one can arrange matters so that $y_{i}=\overline{x_{i}}$ and
$\overline{h_{i}}=-h_{i}$. In this case, the real span of the Cayley transform%
\begin{equation}
\frak{c}_{j} : \left\{  x_{j},h_{j},y_{j}\right\}  \rightarrow
 \left\{ \frac{1}{2}\left( x_{j}+y_{j}-ih_{j}\right) ,-i\left( x_{j}-y_{j}\right) ,
\frac{1}{2}\left( x_{j}+y_{j}+ih_{j}\right) \right\} \quad  \label{cayley} 
\end{equation}
will be a subalgebra $\frak{s}_{i}$ of $\frak{g}_{\mathbb{R}}$ that is isomorphic to
$\frak{s}l\left(  2,\mathbb{R}\right)  $ and moreover%
\[
\left[  \frak{s}_{i},\frak{s}_{j}\right]  =0\qquad,\qquad i\neq j \quad .
\]
Since the semisimple element $h_{i}$ of the original triple $\left\{
x_{i},h_{i},y_{i}\right\}  $ lies in $i\frak{k}_{\mathbb{R}}$ and these all commute we
must have $n\leq rank\left(  K\right)  $. On the other hand, since the
semisimple element $h_{i}^{\prime}$ of the Cayley transform of $\left\{
x_{i},h_{i},y_{i}\right\}  $ is a semisimple element of $\frak{p}_{\mathbb{R}}$, we must
have $n\leq rank\left(  G/K\right)  $. And so it's rather interesting that in
all cases except $SU^{\ast}\left(  2n\right)  $ we're getting the maximal
possible (from this simple argument) number of commuting triples in
$\frak{g}_{\mathbb{R}}$. In the case of $SU^{\ast
}\left(  2n\right)  $, however, the number of $\gamma_{i}$ is $\left[  \frac
{n}{2}\right]  $ , while $rank\left(  K\right)  =n-1$ and $rank\left(
G/K\right)  =n$. We note that there is another circumstance that distinguishes $SU^{\ast
}\left(  2n\right)  $ from the other simple noncompact Lie groups of classical
type: in the case of $SU^{\ast}\left(  2n\right)  $ and only in the case
of $SU^{\ast}\left(  2n\right)  $, there are actually two weights in
$\Delta\left(  \frak{h};\frak{g}\right)  $ that restrict to $\widetilde{\beta
}\in\frak{t}^{\ast}$; that is to say, for $SU^{\ast}\left(  2n\right)  $, and
only $SU^{\ast}\left(  2n\right)  $, there is a pair of complex roots
$\beta,\theta^{\ast}\beta\in\Delta\left(  \frak{h};\frak{g}\right)  $ such
that $\left.  \beta\right|  _{\frak{t}}=\widetilde{\beta}=\left.  \theta
^{\ast}\beta\right|  _{\frak{t}}$. Although, by and large, it is rare
that $\widetilde{\beta}\in\frak{t}^{\ast}$ corresponds to a pair of $\theta
$-conjugate roots in $\Delta\left(  \frak{h};\frak{g}\right)  $ rather than a
single imaginary root, in \S 3 we posit both possibilities on an equal footing.

\section{Families of multiplicity-free $K_{\mathbb{C}}$-orbits}

Let $G$ be a noncompact real semisimple Lie group and let $\Gamma=\left\{
\gamma_{1},\ldots,\gamma_{n}\right\}  $ be a sequence of strongly orthogonal
noncompact weights as constructed in the preceding section. We'll now
associate to $\Gamma$ a corresponding sequence $\left\{  \mathcal{O}%
_{1},\ldots,\mathcal{O}_{n}\right\}  $ of $K_{\mathbb{C}}$-orbits in
$\mathcal{N}_{\frak{p}}$.

We begin by choosing representative elements $x_{i}\in\frak{p}_{\gamma_{i}}$.
As these are nilpotent elements of $\frak{g}$, via a standard construction we
can associate an $S$-triple; that is to say, we can find elements $h_{i}%
,y_{i}\in\frak{g}$ so that for the triple $\left\{  x_{i},h_{i},y_{i}\right\}
$ the following commutation relations are satisfied:
\begin{equation}
\left[  h_{i},x_{i}\right]  =2x_{i},\quad\left[  h_{i},y_{i}\right]
=-2y_{i}\quad,\quad\left[  x_{i},y_{i}\right]  =h_{i}\quad. \label{1}%
\end{equation}
In fact, we can choose $y_{i}\in\frak{p}_{-\gamma_{i}}$ and $h_{i}\in
\frak{t}\subset\frak{k}$ so that $\left\{  x_{i},h_{i},y_{i}\right\}  $ is a
normal triple in $\frak{g}$ ; that is to say, $\left\{  x_{i},h_{i}%
,y_{i}\right\}  $ satisfy both (\ref{1}) and
\begin{equation}
\theta\left(  x_{i}\right)  =-x_{i}\quad,\quad\theta\left(  y_{i}\right)
=-y_{i}\quad,\quad\theta\left(  h_{i}\right)  =h_{i}\quad\label{2}%
\end{equation}
and $\frak{s}_{i}=span_{\mathbb{R}}\left(  x_{i},h_{i},y_{i}\right)  $ is a
$\theta$-stable $\frak{sl}_{2}$-subalgebra of $\frak{g}$ isomorphic to
$\frak{sl}\left(  2,\mathbb{R}\right)  $. Moreover, since the $\gamma_{i}$'s
are strongly orthogonal, the corresponding $\frak{s}_{i}$'s will be mutually
centralizing; i.e., $\left[  \frak{s}_{i},\frak{s}_{j}\right]  =0$ if $i\neq
j$.

We now set%
\begin{align}
X_{i}  &  =x_{1}+x_{2}+\cdots+x_{i} \quad , \nonumber \\
H_{i}  &  =h_{1}+h_{2}+\cdots+h_{i} \quad , \label{XHY} \\
Y_{i}  &  =y_{1}+y_{2}+\cdots+y_{i} \nonumber
\end{align}
and
\[
\mathcal{O}_{i}=K_{\mathbb{C}}\cdot Y_{i}\subset\mathcal{N}_{p}.
\]
$\left\{  X_{i},H_{i},Y_{i}\right\}  $ is easily seen to be another normal
$S$-triple in $\frak{g}$.
\footnote{\emph{Note added in proof.} After submitting
this article, the author learned of a paper by Muller, Rubenthaler and Schiffmann
\cite{mrs} wherein a similar construction of orbits is made.
The situation studied in that paper is where $\frak{g}$ is a complex semisimple
Lie algebra with a $\mathbb{Z}$ grading of the form 
$\frak{g} = \frak{g}_{-1} \oplus \frak{g}_{0} \oplus \frak{g}_{1}$, and the orbits
of the adjoint group of $\frak{g}_{0}$ in $\frak{g}_{1}$. 
As in the Kato-Ochiai \cite{kato-ochiai} paper, however,
sequences of strongly orthogonal weights appear there as an 
auxiliary device arising from an underlying Hermitian-symmetric structure.
By way of contrast, we remark that in the present paper the sequences of strongly 
orthogonal roots are employed as a constructive principle.} 

\bigskip

We'll now show that the orbits $\mathcal{O}_{i}$, $i=1,\ldots,n$ are all multiplicity-free.

To accomplish this, we need to first elaborate a bit more on the setup in
\S 2. Let $G$ be a connected noncompact real semisimple Lie group. Let $K$ be
a maximal compact subgroup, $\theta$ a corresponding Cartan involution and
$\frak{g}=\frak{k}+\frak{p}$, the corresponding Cartan decomposition of the
complexification of the Lie algebra of $G$. Choose a Cartan subalgebra
$\frak{t}$ of $\frak{k}$, and extend it to a $\theta$-stable Cartan subalgebra
$\frak{h}=\frak{t}+\frak{a}$ of $\frak{g}$. Choose a positive system
$\Delta^{+}\left(  \frak{t};\frak{k}\right)  $ for $\Delta\left(
\frak{t};\frak{k}\right)  $ and extend it to a positive system $\Delta
^{+}\left(  \frak{h};\frak{g}\right)  $ of $\Delta\left(  \frak{h}%
;\frak{g}\right)  $ in such a way that
\[
\left.  \alpha\right|  _{\frak{t}}\in\Delta^{+}\left(  \frak{t};\frak{k}%
\right)  \quad\Longrightarrow\quad\alpha\in\Delta^{+}\left(  \frak{h}%
;\frak{g}\right) \quad .
\]
We adopt a Chevalley basis $\left\{  E_{\alpha}\mid\alpha\in\Delta\left(
\frak{h};\frak{g}\right)  \right\}  $ so that
\begin{align}
\left[  E_{\alpha},E_{-\alpha}\right]   &  =H_{\alpha} \quad , \nonumber\\
\left[  E_{\alpha},E_{\gamma}\right]   &  =\left\{
\begin{array}
[c]{cc}%
N_{\alpha,\gamma}E_{\alpha+\gamma} & \text{if }\alpha+\gamma\in\Delta\left(
\frak{h};\frak{g}\right)  \quad , \\
0 & \text{if }\alpha+\gamma\notin\Delta\left(  \frak{h};\frak{g}\right) \quad ,
\end{array}
\right.  \label{3}\\
\left[  H_{\alpha},E_{\gamma}\right]   &  =\left\langle \alpha,\gamma
\right\rangle E_{\gamma}\quad . \nonumber 
\end{align}
and define the induced mapping $\theta^{\ast}:\Delta\left(  \frak{h}%
;\frak{g}\right)  \rightarrow\Delta\left(  \frak{h};\frak{g}\right)  $ and
numbers $\rho_{\alpha}$ by means of the formulas
\begin{align*}
\theta H_{\alpha} &  =H_{\theta^{\ast}\alpha} \quad ,\\
\theta E_{\alpha} &  =\rho_{\alpha}E_{\theta^{\ast}\alpha} \quad .%
\end{align*}
We set
\begin{align*}
\Delta_{0} &  =\left\{  \alpha\in\Delta\mid\theta^{\ast}\alpha=\alpha\right\}
=\text{ the set of pure imaginary roots} \quad , \\
\Delta_{1} &  =\left\{  \alpha\in\Delta\mid\alpha\notin\Delta_{0}\right\}
=\text{ the set of complex roots} \quad , %
\end{align*}
and
\[
\Delta_{0,\pm}=\left\{  \alpha\in\Delta_{0}\mid\theta E_{\alpha}=\pm
E_{\alpha}\right\} \quad .
\]
$\Delta_{0},_{+}$ and $\Delta_{0,-}$ are, respectively, the sets of,
compact imaginary roots and noncompact imaginary roots. We
remark that since we have set up $\frak{h}$ as a maximally compact Cartan
subalgebra, there are no real roots in $\Delta\left(  \frak{h};\frak{g}%
\right)  $.

Let $\beta\in\Delta\left(  \frak{h};\frak{g}\right)  $ be a root such that
\begin{equation}
\left(  1-\theta\right)  \frak{g}_{\beta}=\frak{p}_{\widetilde{\beta}}%
\equiv\text{ the highest weight space of the representation of }K\text{ on
}\frak{p}\quad. \tag{*}%
\end{equation}
There are several situations that we shall reduce to two basic cases:

\begin{itemize}
\item [ (i) ]$\beta$ is a complex root;

\item[ (ii) ] $\beta$ is a pure imaginary noncompact root.
\end{itemize}

When $rank\left(  \frak{g}\right)  =rank\left(  \frak{k}\right)
$, all roots in $\Delta\left(  \frak{h};\frak{g}\right)  =\Delta\left(
\frak{t;g}\right)  $ will be pure imaginary and so $\beta$ will be a
non-compact imaginary root. When $rank\left(  \frak{g}\right)  >rank\left(  \frak{k}\right)
$, $\beta$ will either be a non-compact pure imaginary root or there will be 
a  $\theta^{\ast}$-conjugate pair
of complex roots $\left\{  \beta,\theta^{\ast}\beta\right\}  $ sharing the
property (*). In the latter case, even though we make an initial choice for
$\beta$, subsequent developments will be manifestly independent of that choice.

\begin{Lemma}
Suppose $\beta$ is a complex root and $\widetilde{\beta}=\left.
\beta\right|  _{\frak{t}}$ is the highest weight of the representation of $K$
on $\frak{p}$. Then both%
\[
\left[  \theta E_{\pm\beta},E_{\mp\beta}\right]  =0
\]
and%
\[
\left[  \theta E_{\pm\beta},E_{\pm\beta}\right]  =0 \quad .
\]
\end{Lemma}

\emph{Proof.} To prove the first, we note that%
\[
\left[  \theta E_{\pm\beta},E_{\mp\beta}\right]  \in\frak{g}_{\pm\left(
\theta^{\ast}\beta-\beta\right)  } \quad .%
\]
However, $\pm\left(  \theta^{\ast}\beta-\beta\right)  $ will be real root, but
for our choice of $\frak{h}$ there are no real roots, and so $\theta
E_{\pm\beta}$ and $E_{\mp\beta}$ must commute.

To prove the second relation, we set%
\begin{align*}
k_{\beta}  &  =\left(  1+\theta\right)  E_{\beta}\in\frak{k}_{\widetilde
{\beta}} \quad , \\
p_{\beta}  &  =\left(  1-\theta\right)  E_{\beta}\in\frak{p}_{\widetilde
{\beta}} \quad .%
\end{align*}
Since $k_{\beta}$ is positive root vector for $K$ and $p_{\beta}$ is the
highest weight of the representation of $K$ on $\frak{p}$ we must have%
\begin{align*}
0  &  =\left[  k_{\beta},p_{\beta}\right]  =\left[  \left(  1+\theta\right)
E_{\beta},\left(  1-\theta\right)  E_{\beta}\right] \\
&  =\left[  E_{\beta},E_{\beta}\right]  +\left[  \theta E_{\beta},E_{\beta
}\right]  -\left[  E_{\beta},\theta E_{\beta}\right]  -\left[  \theta
E_{\beta},\theta E_{\beta}\right] \\
&  =2\left[  \theta E_{\beta},E_{\beta}\right] \quad .
\end{align*}
Similarly, $\left[  k_{-\beta},p_{-\beta}\right]  $ implies $\left[  \theta
E_{-\beta},E_{-\beta}\right]  =0.$

\qed

\begin{Lemma}
Suppose $\beta\in\Delta\left(  \frak{h};\frak{g}\right)  $ is such
that $\widetilde{\beta}=\left.  \beta\right|  _{\frak{t}}$ is a highest
weight of the representation of $K$ on $\frak{p}$. Let $x_{\beta},h_{\beta
},y_{\beta}$ be defined by%
\begin{align*}
x_{\beta} &  =E_{\beta} \quad , \\
h_{\beta} &  =\frac{2}{\left\langle \beta,\beta\right\rangle }H_{\beta} \quad ,\\
y_{\beta} &  =E_{-\beta} \quad ,%
\end{align*}
if $\beta$ is pure imaginary, or
\begin{align*}
x_{\beta} &  =\left(  1-\theta\right)  E_{\beta} \quad , \\
h_{\beta} &  =\frac{2}{\left\langle \beta,\beta\right\rangle }\left(
1+\theta\right)  H_{\beta} \quad , \\
y_{\beta} &  =\frac{2}{\left\langle \beta,\beta\right\rangle }\left(
1-\theta\right)  E_{\beta}%
\end{align*}
if $\beta$ is complex. Then $\left\{  x_{\beta},h_{\beta},y_{\beta}\right\}  $
is a normal $S$-triple in $\frak{g}$ and
\begin{equation}
z\in\frak{p}\text{ and }\left[  h_{\beta},z\right]  =-2z\quad\Longrightarrow
\quad z\in span_{\mathbb{C}}\left(  y_{\beta}\right)  \quad .\tag{**}%
\end{equation}
\end{Lemma}

\emph{Proof.} \ An obvious calculation using the commutation relations
(\ref{3}), and Lemma 3.1 in the case when $\beta$ is complex, confirms that
\[
\left[  h_{\beta},x_{\beta}\right]  =2x_{\beta}\quad,\quad\left[  h_{\beta
},y_{\beta}\right]  =-2y_{\beta}\quad,\quad\left[  x_{\beta},y_{\beta}\right]
=h_{\beta}%
\] and so $\left\{  x_{\beta},h_{\beta},y_{\beta}\right\}  $ is an $S$-triple. It
is also obvious that
\[
\theta x_{\beta}=-x_{\beta}\quad,\quad\theta y_{\beta}=-y_{\beta}\quad
,\quad\theta h_{\beta}=h_{\beta}%
\]
and so $x_{\beta},y_{\beta}\in\frak{p}$, $h_{\beta}\in\frak{k}$; hence
$\left\{  x_{\beta},h_{\beta},y_{\beta}\right\}  $ is a normal $S$-triple. Of
course, $x_{\beta}$ and $y_{\beta}$ live, respectively, in the $+2$ and $-2$
eigenspaces of $h_{\beta}$.

To prove (**) we have to show that no other weight vector in $\frak{p}$ can
live in the $-2$-eigenspace of $\frak{h}_{\beta}$. We shall handle the cases
$\beta$ is a complex root or a noncompact imaginary root separately.

\emph{Case (i)}. Assume $\beta$ is a non-compact complex root and set
\[
k_{\beta}=\left(  1+\theta\right)  E_{\beta}\qquad,\qquad k_{-\beta}=\frac
{2}{\left\langle \beta,\beta\right\rangle }\left(  1+\theta\right)  E_{-\beta} \quad .%
\]
Then it is easy to check that $\left\{  k_{\beta},h_{\beta},k_{-\beta
}\right\}  $ is a $\theta$-stable $S$-triple in $\frak{g}$ with the same
semisimple element as that of $\left\{  x_{\beta},h_{\beta},y_{\beta}\right\}
$. Suppose $\widetilde{\alpha}\in\frak{t}^{\ast}$ is a weight of the
representation of $\frak{k}$ on $\frak{p}$, and $\widetilde{\alpha}\neq
\pm\widetilde{\beta}=\pm\left.  \beta\right|  _{\frak{k}}$. \ Having chosen a
positive system for $\Delta\left(  \frak{h};\frak{g}\right)  $ subordinate to
that of $\Delta\left(  \frak{t};\frak{k}\right)  $, we can regard
$\widetilde{\alpha}$ as, respectively, a ``positive'' or ``negative'' weight
of $\frak{p}$, depending on whether $\widetilde{\alpha}=\left.  \alpha\right|
_{\frak{t}}$ for some $\alpha\in\Delta^{\pm}\left(  \frak{h};\frak{g}\right)
$. Assume $\widetilde{\alpha}$ is ``positive''; then $\left[  k_{\beta
},z\right]  \in\frak{p}_{\widetilde{\alpha}+\widetilde{\beta}}=\left\{
0\right\}  $ since $\widetilde{\beta}$ is the highest weight of the
representation of $K$ on $\frak{p}$, and $\left[  k_{-\beta},\left[
k_{-\beta},z\right]  \right]  \in\frak{p}_{\widetilde{\alpha}-2\widetilde
{\beta}}=\left\{  0\right\}  $ since $-\beta$ is the lowest weight of
$\frak{p}$. A similar (albeit up-side-down) argument shows that if
$\widetilde{\alpha}$ is ``negative'', then $\left[  k_{-\beta},z\right]
=0=\left[  k_{\beta},\left[  k_{\beta},z\right]  \right]  $. And so, in either
case, the representation of $\frak{sl}\left(  2,\mathbb{R}\right)  $ generated by the
action of $k_{\pm\beta}$ on $z\in\frak{p}_{\widetilde{\alpha}}$ is at most
$2$-dimensional, and so the lowest possible eigenvalue of $h_{\beta}$ is $-1$.

\emph{Case (ii)}. Assume $\beta$ is a pure imaginary noncompact root. In this
case, the $S$-triple $\left\{  x_{\beta},h_{\beta},y_{\beta}\right\}  $ is
just a renormalization of the Weyl triple $\left\{  E_{\beta},H_{\beta
},E_{-\beta}\right\}  $. Since each weight $\widetilde{\alpha}$ of $\frak{p}$
either comes from a pair of $\theta^{\ast}$-conjugate complex roots
$\alpha,\theta^{\ast}\alpha$ $\in\Delta_{1}\left(  \frak{h};\frak{g}\right)  $
or corresponds more or less directly to a unique pure imaginary non-compact
root, it will suffice to show that for any root $\alpha\in\Delta\left(
\frak{h};\frak{g}\right)  $ such that $\alpha\neq\pm\beta$, the maximal length
of a $\beta$-string through $\alpha$ is $2$.

Suppose $\alpha$ is root in $\Delta_{1}^{+}\left(  \frak{h};\frak{g}\right)  $
then $\left[  E_{\alpha},E_{\beta}\right]  =0$, since otherwise $\alpha+\beta$
would be a complex root and we'd have a root vector with a non-zero projection
to $\frak{p}$ and a weight higher than $\beta$. Therefore, for any root
$\alpha\in\Delta_{1}^{+}\left(  \frak{h};\frak{g}\right)  $, $\alpha$ would
have to be at the top of a (perhaps trivial) $\beta$-string. And so the
situation we have to worry about when $\alpha\in\Delta_{1}^{+}\left(
\frak{h};\frak{g}\right)  $ is when the string is $\alpha,\alpha-\beta
,\alpha-2\beta$ or longer. If $\alpha-2\beta$ is a root, it must be a complex
root and so $\left.  \left(  \alpha-2\beta\right)  \right|  _{\frak{t}}$
\ must be a $\frak{t}$-weight of $\frak{p}$. But $\left.  \left(
\alpha-2\beta\right)  \right|  _{\frak{t}}$ \ is a weight lower than $-\beta$,
the lowest weight of $\frak{p}$. Hence, we have a contradiction if
$\alpha-2\beta\in\Delta\left(  \frak{h},\frak{g}\right)  $. If $\alpha
\in\Delta_{1}^{-}\left(  \frak{h};\frak{g}\right)  $, an analogous (albeit
upside-down)\ arguments show that neither $\alpha-\beta$ or $\alpha+2\beta$
can be roots in $\Delta\left(  \frak{h};\frak{g}\right)  $. \ 

Our last concern then would be the possible existence of a $\beta$-string
$\ldots,\alpha-\beta,\alpha,\alpha+\beta,\ldots$ through a non-compact pure
imaginary root $\alpha\neq-\beta$ such that $\left[  h_{\beta},E_{\alpha
}\right]  =-2E_{\alpha}$. In this case we'd have $\left[  x_{\beta},\left[  x_{\beta
},E_{\alpha}\right]  \right]  \in\frak{p}_{2\beta+\alpha}$, which is
impossible since $\beta$ is the highest weight of $\frak{p}$.

\qed

\begin{Corollary}
If $z\in\frak{p}$ and $\left[  H_{i},z\right]  =-2z$, then $z\in
span_{\mathbb{C}}\left(  y_{1},\ldots,y_{i}\right)  $.
\end{Corollary}

\emph{Proof. }Since the weights $\gamma_{i}$ are all in the Weyl group orbit
of the highest weight $\widetilde{\beta}$ of $\frak{p}$, for each $i$ we can
choose a positive systems so that $x_{i}$ is a highest weight vector in
$\frak{p}$. It then follows from the preceding lemma, that for each
$i=1,\ldots,n$, the lowest eigenvalue of $h_{i}$ will be $-2$ and $\left[
h_{i},z\right]  =-2z$ will imply that $z\in\mathbb{C}y_{i}$. Since the $h_{i}$
are simultaneously diagonalizable, we can conclude that the smallest
eigenvalue of $H_{i}=h_{1}+\cdots+h_{i}$ will be $-2$ and that if $\left[
H_{i},z\right]  =-2z$ then we must have $z\in span_{\mathbb{C}}\left(
y_{1},\ldots,y_{i}\right)  $. \qed

Let $\frak{t}_{1}=span_{\mathbb{C}}\left(  h_{1},\ldots,h_{i}\right)  $ and
let $\frak{t}_{0}$ be the orthogonal complement of $\frak{t}_{1}$ in
$\frak{t}$ (with respect to the Killing form). Let $\frak{m}_{i}$ be the
subalgebra of $\frak{k}$ generated by root spaces $\frak{k}_{\alpha}$ such
that $\left.  \alpha\right|  _{\frak{t}_{1}}=0$.

\begin{Lemma} Let $\overline{\frak{n}}_{i}$ be the
direct sum of the negative eigenspaces of $ad\left(  H_{i}\right)  $ in
$\frak{k}$. Then
$\frak{m}_{i}+ \overline{\frak{n}}_{i} \subseteq\frak{k}^{Y_{i}}$.
\end{Lemma}

\emph{Proof}. Since $-2$ is the lowest eigenvalue of $H_{i}$, certainly 
$\overline{\frak{n}}_{i} \subset \frak{k}^{Y_{i}}$.
Suppose $k_{\alpha}$ is a root vector corresponding to a root
space $\frak{t}_{\alpha}\subset\frak{m}_{i}$. \ We then have $\left[
H_{i},k_{\alpha}\right]  =0$ and so $k_{\alpha}$ will preserve the $\left(
-2\right)  $-eigenspace of $H_{i}$. \ By the preceding lemma
\[
\left(  -2\right)  \text{-eigenspace of }H_{i}=span_{\mathbb{C}}\left(
y_{1},\ldots,y_{i}\right) \quad .
\]
Because the $y_{j}$, $j=1,\ldots,i$ are weight vectors corresponding to
multiplicity-free weights of $\frak{p}$, we must have%
\[
\left[  k_{\alpha},y_{j}\right]  =cy_{k}\qquad\text{for some \thinspace}%
k\in\left\{  1,\ldots,i\right\}  -j\text{ and some }c\in\mathbb{C}\quad.
\]
But if $c\neq0$, $\gamma_{k}$ will not be strongly orthogonal to $\gamma_{j}$;
for otherwise we would have $\gamma_{k}-\gamma_{j}=\alpha\in\Delta\left(
\frak{t};\frak{k}\right)  $. We conclude that $Y_{i}$ commutes with every root
vector in $\frak{m}_{i}$ and so, since $\frak{m}_{i}$ is semisimple, $\left[
Y_{i},\frak{m}_{i}\right]  =0$.  \qed

It now follows readily from results of Servedio \cite{servedio} and Kimel'fel'd-Vinberg
\cite{kimelfeld-vinberg} that the $K_{\mathbb{C}}$-orbit through $Y_{i}$ is multiplicity
free. \footnote{A little more explicitly, the argument would proceed as
follows. It is easy to see that the tangent space to $Y_{i}$ is generated by a
parabolic subalgebra of $\frak{k}$ corresponding to the non-negative
eigenspaces of $ad_{\frak{k}}H_{i}$. The preceding lemma then allows one to
whittle this parabolic down to a Borel subalgebra. One then applies
\par
\textbf{Theorem} (\cite{servedio}, \cite{kimelfeld-vinberg}) 
Let $K$ be a connected reductive algebraic group
acting on an irreducible affine variety $M$. Then $\mathbb{C}\left[  M\right]
$ is multiplicity-free if and only if there exists a Borel subgroup $B\subset K$
admitting an open orbit in $M$.}
However, we shall instead apply
algebraic Frobenius reciprocity; so that we can not only demonstrate that
$\mathbb{C}\left[  \overline{\mathcal{O}_{i}}\right]  $ is multiplicity-free,
but we can also identify the $K$-types.

\begin{Theorem}
[Kostant, \cite{kostant}]Suppose $G$ is a reductive algebraic group and $V$ is an
irreducible $G$-module. For $\ x\in V,$ let $G^{x}$ denote the stabilizer of
$x$ in $G$. \ If we denote by $O_{x}$ the $G$-orbit through $x$, and by
$\mathbb{C}\lbrack \mathcal{O}_{x}\rbrack  $ the ring of everywhere-defined
rational functions on $\mathcal{O}_{x}$, by $V_{\lambda}$, $\lambda
\in\widehat{G}$, the representation space of an irreducible
finite-dimensional representation of $G$ and by 
$\widetilde{V_{\lambda}}$ the dual module of $V_{\lambda}$, then 
\[
\text{multiplicity of }\lambda\text{ in }\mathbb{C}\lbrack  \mathcal{O}%
_{x}\rbrack   = \dim\widetilde{V_{\lambda}}^{G^{x}}%
\]
where 
$\widetilde{V_{\lambda}}^{G^{x}}$ is the space of vectors in $\widetilde
{V_{\lambda}}$ that are fixed by $G^{x}$.
\end{Theorem}

\textbf{Remark:} In the preceding theorem the stabilizer $G^{x}$ of $x$ need
not be reductive.

\begin{Corollary}
$\mathbb{C}\lbrack  \mathcal{O}_{i}\rbrack  $ is multiplicity-free and
$V_{\lambda}$ is a $K$-type in $\mathbb{C}\lbrack  \mathcal{O}_{i}\rbrack  $
then $\lambda\in span_{\mathbb{R}}\left(  \gamma_{1},\ldots,\gamma_{i}\right)
$.
\end{Corollary}

\emph{Proof.} By Lemma 3.4 the stabilizer of $Y_{i}$ in $\frak{k}$ contains
$\frak{m}_{i}+\overline{\frak{n}}_{i}$. It is easy to see that the semisimple element ${h_j}$ in the normal
$S$-triple $\left\{x_{j},h_{j},y_{j}\right\}$ is an element of
$\frak{t}$ such that $\lbrack h_{j} , k_{\alpha} \rbrack = 
\langle \alpha , \gamma_{j} \rangle k_{\alpha} $ for any root vector $k_{\alpha}
\in \frak{k}_{\alpha}$, $\alpha \in \Delta \left(\frak{t};\frak{k}\right)$.  
Since the $\gamma_{i}$ are chosen such that each $\frak{t}$-weight $\gamma_{1}+\cdots+\gamma_{i}$ is
 dominant and since $H_{i} \equiv h_{1} + \cdots + h_{i}$, it follows that all the negative root vectors of
$\frak{k}$ will be contained in $\frak{m}_{i}+\overline{\frak{n}_{i}}$. Hence,
a nonzero element of $\widetilde{V_{\lambda}}^{K^{Y_{i}}}$ will be a lowest
weight vector that is also $\frak{m_{i}}$-invariant. Algebraic Frobenius
reciprocity then implies that if a $K$-type $V_{\lambda}$ appears in
$\mathbb{C}\lbrack  \mathcal{O}_{i}\rbrack  $ then the lowest weight vector of
$\widetilde{V_{\lambda}}$ must be $\frak{m}_{i}$-invariant. This in turn
implies that $-\lambda$ (and so $\lambda$) is not supported on $\frak{t}_{0}$.
Thus, we must have%
\[
\lambda=a_{1}\gamma_{1}+\cdots+a_{i}\gamma_{i}\in\frak{t}_{1}^{\ast} \quad .%
\]
And, of course, since the space of lowest weight vectors in $\widetilde
{V_{\lambda}}$ will be $1$-dimensional, algebraic Frobenius reciprocity also
tells us that $\mathbb{C}\lbrack  \mathcal{O}_{i}\rbrack  $ is
multiplicity-free. \qed

\begin{Proposition}
Let $n=\left| \Gamma \right|$. (i) If $i < n$, then $V_{\lambda }$ is a $K$-type in $%
\mathbb{C}\lbrack \mathcal{O}_{i}\rbrack $, if and only if its highest
weight is of the form
\begin{equation*}
\lambda =a_{1}\gamma _{1}+\cdots +a_{i}\gamma _{i}\qquad 
\end{equation*}
with $a_{j}\in \mathbb{N}$ and $a_{1}\geq a_{2}\geq \cdots \geq a_{i}\geq 0$%
. (ii) $V_{\lambda }$ is a $K$-type in $\mathbb{C}\lbrack \mathcal{O}%
_{n}\rbrack $, if and only if $\lambda $ is of the form
\begin{equation*}
\lambda =a_{1}\gamma _{1}+\cdots +a_{n}\gamma _{n}
\end{equation*}
with $a_{j}\in \mathbb{Z}$ and
\begin{eqnarray*}
a_{1} &\geq &a_{2}\geq \cdots \geq a_{n-1}\geq a_{n}\quad \quad \quad ,\quad 
\text{if }\Sigma =\left( a_{n}\right) ^{m_{a}}\ \text{\ or }\ \left(
A_{n}\right) ^{m_{A}}\quad ; \\
a_{1} &\geq &a_{2}\geq \cdots \geq a_{n-1}\geq \left| a_{n}\right| \geq
0\quad ,\quad \text{if }\Sigma =\text{ }\left( d_{n}\right) ^{m_{d}}\quad ;
\\
a_{1} &\geq &a_{2}\geq \cdots \geq a_{n-1}\geq a_{n}\geq 0\quad ,\quad \text{%
otherwise}\quad .
\end{eqnarray*}
\end{Proposition}

\emph{Proof}. From the Corollary above, we know that if $V_{\lambda }$ is a $%
K$-type in $\mathbb{C}\lbrack \mathcal{O}_{i}\rbrack $ then its highest
weight must be of the form 
\begin{equation*}
\lambda =a_{1}\gamma _{1}+\cdots +a_{i}\gamma _{i}
\end{equation*}
We first show that coefficients $a_{k}$ must be integers. Note that 
\begin{equation*}
\exp \left( i\pi h_{j}\right) \cdot Y_{i}=\exp \left( -2i\pi \right)
Y_{i}=Y_{i}\qquad \text{for }1\leq j\leq i
\end{equation*}
and so $k_{j}\equiv \exp \left( i\pi h_{j}\right) \in K^{Y_{i}}$. On the
other hand, $\lambda =a_{1}\gamma _{1}+\cdots +a_{i}\gamma _{i}$ and $%
v_{-\lambda }$ is the lowest weight vector of $\widetilde{V_{\lambda }}$ we
have 
\begin{equation*}
k_{j}\cdot v_{-\lambda }=\exp \left( -2i\pi a_{j}\right) v_{-\lambda }\text{
.}
\end{equation*}
Thus, the lowest weight vector of $\widetilde{V_{\lambda }}$ will not be
stabilized by $k_{i}\in K^{Y_{i}}$ unless $a_{j}\in \mathbb{Z}$ for $%
j=1,\ldots ,i$.

We now to prove the necessity of the ordering of the coefficients $a_{i}$. This is just a
consequence of requirement that the highest weight $\lambda $ be dominant.
Since, in all cases, for all $1\leq j<k\leq n$,  
\begin{equation*}
\gamma _{j}-\gamma _{k}\text{ or }\frac{1}{2}\gamma _{i}-\frac{1}{2}\gamma
_{j}\in \Sigma ^{+}
\end{equation*}
The dominance condition on $\lambda =a_{1}\gamma _{1}+\cdots +a_{i}\gamma
_{i}$, leads to
\begin{equation*}
a_{j}\geq a_{j+1}\qquad i=j=1,\ldots ,n-1
\end{equation*}
If the restricted root system $\Sigma $ contains a $d_{n}$ factor, then we
must have in addition
\begin{equation*}
\left\langle \lambda ,\gamma _{j}+\gamma _{k}\right\rangle \geq 0
\end{equation*}
which, together with $\left\langle \lambda ,\gamma _{n-1}-\gamma
_{n}\right\rangle \geq 0$,  implies in particular that 
\begin{equation*}
a_{n-1}\geq \left| a_{n}\right| \quad .
\end{equation*}
In all other cases, $\Sigma $ contains either a $b_{n}$ or $C_{n}$ factor,
and this leads to the requirement that
\begin{equation*}
\left\langle \lambda ,\gamma _{i}\right\rangle \geq 0\qquad \text{for all }%
i=1,\ldots ,n
\end{equation*}
which forces all the coefficients $a_{i}$ to be non-negative. 

At this point we have seen that if $\lambda$ is the highest weight
of a $K$-type in $\mathbb{C}\left[ \mathcal{O}_{i} \right]$, then
(1) $\lambda$ must lie in the span of the $\gamma_{j}$, $j\le i$, (2)
the coefficients $a_{1}, \ldots , a_{i}$ of $\lambda$ with respect to 
$\gamma_{1}, \ldots , \gamma_{i}$
must satisfy certain integrality conditions so that each 
$\exp\left(2\pi h_{j} \right) \in K^{Y_{i}}$ acts trivially on the corresponding
highest weight vector and (3) the coefficients must be ordered in such a way
that $\lambda$ is dominant. What is not yet clear is that these restrictions
on $\lambda$ are sufficient to place $V_{\lambda}$ 
in $\mathbb{C}\left[ \mathcal{O}_{i} \right]$. However, it is easy to see that
the integrality conditions we have imposed on $\lambda$ are actually stronger than
those needed to guarantee that $\lambda$ is a weight of a representation of $K$.
Thus, we have enumerated all possible finite dimensional representations of $K$
with a $K^{Y_i}$-fixed vector. Sufficiency now follows from Theorem
3.5. 
 \qed

The preceding proposition tells us exactly which $K$-types occur in the
ring $\mathbb{C}\lbrack \mathcal{O}_{i}\rbrack $. However, we are actually most 
interested in the ring of regular functions $\mathbb{C}\left[ \overline{\mathcal{O}_{i}}%
\right] $ on the closure of the orbit. Clearly,
\begin{equation*}
\mathbb{C}\left[ \overline{\mathcal{O}_{i}}\right] \subset \mathbb{C}\lbrack 
\mathcal{O}_{i}\rbrack \quad .
\end{equation*}
We will now show that each $K$-type $V_{\lambda }$ occuring in $\mathbb{C}%
\lbrack \mathcal{O}_{i}\rbrack $ also occurs in $\mathbb{C}\left[ \overline{%
\mathcal{O}_{i}}\right] $.

\begin{Theorem}
[Kumar, \cite{kumar}]Let $\frak{g}$ be a finite-dimensional semisimple Lie algebra
and let $\frak{h}$ be a Cartan subalgebra of $\frak{g}$. Let $V_{\lambda}$ denote the
irreducible finite-dimensional representation of $\frak{g}$ with highest weight
$\lambda$. For any weight $\lambda \in \frak{h}^{\ast}$, let 
$\overline{\lambda}$ denote the unique dominant weight in the Weyl group orbit of $\lambda$.  
Then, for any pair
$\lambda,\mu$ of dominant weights and any $w$ in the Weyl group of $\frak{g}$,
the irreducible $\frak{g}$-module $V_{\overline{\lambda+\omega\mu}}$ 
occurs with multiplicity exactly one in the $\frak{g}%
$-submodule $U\left(  \frak{g}\right)  \cdot\left(  e_{\lambda}\otimes
e_{w\mu}\right)  $ of $V_{\lambda}\otimes V_{\mu}$; where $e_{\lambda}$ and
$e_{w\mu}$ are, respectively, weight vectors in the $\lambda$-weight space of
$V_{\lambda}$ and the $w\mu$-weight space of $V_{\mu}$).
\end{Theorem}

\begin{Remark}
The statement of the theorem is known as Kostant's strengthened
Parthasarathy-Ranga Rao-Varadarajan conjecture.
\end{Remark}

\begin{Lemma}
Let $\omega_{j}=\gamma_{1}+\gamma_{2}+\cdots+\gamma_{j}$, $1\leq j\leq i$.
Then the $K$-type $V_{\omega_{j}}$ occurs in $S^{j}\left(  \frak{p}\right)  $
and the monomial $x_{1}\cdots x_{j}\in S^{j}\left(  \frak{p}\right)  $ has a
non-trivial projection onto the highest weight space of $V_{\omega_{j}}$
\end{Lemma}

\emph{Proof.} The case when $j=1$ is trivial, since $x_{i}$ is a highest
weight vector of $\frak{p}$. We now proceed by induction on $i$. And all that
this
requires is the preceding theorem with the identification of $e_{\lambda}$
with the projection of $x_{1}\cdots x_{j}$ onto the highest weight vector of
$V_{\omega_{j}}$ (the inductive hypothesis) and the identification of
$e_{w\mu}$ with $x_{j+1}$ which, by our construction, is always 
extremal weight vector of $\frak{p}$. \qed

\begin{Theorem}
The $K$-type decomposition of $\mathbb{C}\left[ \overline{\mathcal{O}_{i}}%
\right] $ is exactly
\begin{equation*}
\mathbb{C}\left[ \overline{\mathcal{O}_{i}}\right] =\bigoplus_{\lambda \in
\Lambda _{i}}V_{\lambda }
\end{equation*}
where
\begin{equation*}
\Lambda _{i}=\left\{ \lambda =a_{1}\gamma _{1}+\cdots +a_{i}\gamma _{i}\mid
a_{j}\in \mathbb{N}\quad ,\quad a_{1}\geq a_{2}\geq \cdots \geq a_{i}\geq
0\right\} 
\end{equation*}
if $i<n$ or $\Sigma \neq (d_{n})^{m_{d}}$. If $i=n$ and $\Sigma
=(d_{n})^{m_{a}}$, then 
\begin{equation*}
\Lambda =\left\{ \lambda =a_{1}\gamma _{1}+\cdots +a_{n}\gamma _{n}\mid
a_{j}\in \mathbb{Z}\quad ,\quad a_{1}\geq a_{2}\geq \cdots \geq a_{n-1}\geq
|a_{n}|\geq 0\right\} \quad .
\end{equation*}
\end{Theorem}

\emph{Proof.} We have already seen that $\lambda \in \Lambda _{i}$ is a
necessary condition for a $K$-type to be in $\mathbb{C}\left[ \mathcal{O}_{i}%
\right] \supset \mathbb{C}\left[ \overline{\mathcal{O}_{i}}\right] $. What
we must prohibit, or otherwise take into account, is the existence of a
rational function that is defined everywhere on $\mathcal{O}_{i}$ but does
not extend to the boundary of $\mathcal{O}_{i}$; in particular, rational
functions that are not defined at $0$. 

On the other hand, if $\phi $ is a polynomial function on $\frak{p}$ that is
supported at some point on $\mathcal{O}_{i}$ then $0\neq \left. \phi \right|
_{\overline{\mathcal{O}_{i}}}\in \mathbb{C}\left[ \overline{\mathcal{O}_{i}}%
\right] $. We now note that the monomials $x_{1}\cdots x_{j}$ are supported
at $Y_{i}$; for
\begin{align*}
\left( x_{1}\cdots x_{j}\right) \left( Y_{i}\right) & =\left\langle
x_{1},Y_{i}\right\rangle \cdots \left\langle x_{j},Y_{i}\right\rangle
=\left\langle x_{1},y_{1}+\cdots +y_{i}\right\rangle \cdots \left\langle
x_{j},y_{1}+\cdots +y_{i}\right\rangle  \\
& =\left\langle x_{1},y_{1}\right\rangle \cdots \left\langle
x_{j},y_{j}\right\rangle  \\
& =1\quad .
\end{align*}
In fact, we can choose an orthogonal basis for $\frak{p}$ such that $%
x_{1},\ldots ,x_{j}$, $j\leq i$, are the only coordinate functions supported
at $Y_{i}$. By the preceding lemma, for each $j=1,\ldots ,i$, there exists a
(homogeneous) highest weight vector $\phi _{\omega _{j}}$ of $V_{\omega _{j}}\subset S^{j}(%
\frak{p})$ of the form
\begin{equation*}
\phi _{\omega _{j}}=x_{1}\cdots x_{j}\ +\ \left( \text{other terms with at
least one factor not among }x_{1},\ldots ,x_{i}\right) \quad .
\end{equation*}
(Note that by the homogeneity of $\phi_{\omega_{i}}$ and the linear independence of the weights $\gamma _{1},\ldots
,\gamma _{i}$, we cannot have a factor $x_{k}$, $j<k\leq i$ occurring in one
of the ``other terms'' without an accompanying coordinate outside of $%
\left\{ x_{1},\ldots ,x_{i}\right\} $.) The ``other terms'' will thus die
upon evaluation at $Y_{i}$ and so we must have
\begin{equation*}
\phi _{\omega _{j}}\left( Y_{i}\right) =1\qquad \text{for}\quad 1\leq j\leq
i\quad .
\end{equation*}
But once we have the highest weight vectors of each $V_{\omega _{j}}$
supported at $Y_{i}$ it is trivial to show that the products of these
highest weight vectors will remain highest weight vectors and continue to be
supported at $Y_{i}\in \mathcal{O}_{i}$. Thus, each of the $K$-types $%
V_{\lambda }$ with
\begin{equation*}
\lambda \in \Lambda _{i}^{\prime }\equiv \left\{ \alpha _{1}\omega
_{1}+\cdots +\alpha _{i}\omega _{i}\mid \alpha _{i}\in \mathbb{N}\right\} 
\end{equation*}
will be supported at $Y_{i}$. We now observe $\Lambda _{i}=\Lambda
_{i}^{\prime }$ So each of the $K$-types specified in the statement of the
theorem definitely appears in $\mathbb{C}\left[ \overline{\mathcal{O}_{i}}%
\right] $. Comparing with Proposition 3.7, we can conclude that in fact so
long as $i\neq n$ and  $\Sigma \neq A_{n}$ or $a_{n}$, we have
\begin{equation*}
\mathbb{C}\left[ \overline{\mathcal{O}_{i}}\right] =\mathbb{C}\left[ 
\mathcal{O}_{i}\right] =\bigoplus_{\lambda \in \Lambda _{i}}V_{\lambda }
\end{equation*}

However, when $i=n$ and  $\Sigma =A_{n}$ or $a_{n}$ we seem to have $K$%
-types $V_{\alpha _{1}\gamma _{1}+\cdots +a_{n}\gamma _{n}}$ with $a_{1}\geq
a_{2}\geq \cdots \geq a_{n}$ and $a_{n}<0$ appearing in $\mathbb{C}\left[ 
\mathcal{O}_{n}\right] $ that have not been accounted for in $\mathbb{C}%
\left[ \overline{\mathcal{O}_{n}}\right] $. In fact, they do not occur in $%
\mathbb{C}\left[ \overline{\mathcal{O}_{n}}\right] ;$ yet they can
nevertheless be easily taken into account. 

Note that the restricted root
system $\Sigma $ is of type $A_{n}$ or $a_{n}$ only when $G/K$ is Hermitian
symmetric (see also \cite{sahi0}). When $Sigma$ is of
this form it is easy to see that $\omega _{n}=\gamma _{1}+\cdots +\gamma _{n}\in 
\frak{t}^{\ast }$ is perpendicular to every root in $\Delta \left( \frak{t};%
\frak{k}\right) $, being supported only on the center of $\frak{k}$. The
corresponding highest weight vector $\phi _{\omega _{n}}\in S^{n}\left( 
\frak{p}\right) $ thus corresponds to a one-dimensional $K$-type which, as
shown above, does not vanish at $Y_{n}.$  Put another way, $\phi _{\omega
_{n}}$ corresponds to a $K$-semi-invariant polynomial that does not vanish
at $Y_{n}\in \mathcal{O}_{n}$. But then it vanishes nowhere on $\mathcal{O}%
_{n}$ (otherwise $\mathcal{O}_{n}$ would have a proper $K$-invariant
subset). Yet being a homogeneous polynomial of degree $n$, it certainly
vanishes at $0\in \overline{\mathcal{O}_{n}}$. Now observe that if $\lambda
=a_{1}\gamma _{1}+\cdots +a_{n}\gamma _{n}\in \Lambda _{n}$, $a_{n}\geq 0$,
and  $\phi _{\lambda }$ is a polynomial corresponding the highest weight
vector of $V_{\lambda }\subset \mathbb{C}\left[ \overline{\mathcal{O}_{i}}%
\right] $ then 
\begin{equation*}
\psi =\left( \phi _{\omega _{n}}\right) ^{-s}\left( \phi _{\lambda }\right) 
\end{equation*}
will be a rational function of highest weight $\left( a_{1}-s\right) \gamma
_{1}+\cdots +\left( a_{n}-s\right) \gamma _{n}$ that is everywhere defined
on $\mathcal{O}_{n}$ but undefined at $0$ if $a_{n}-s<0$.  Evidently, such
functions account for the highest weight vectors of all the $K$-types in $%
\mathbb{C}\left[ \mathcal{O}_{n}\right] $ that have not already been shown
to appear in $\mathbb{C}\left[ \overline{\mathcal{O}_{n}}\right] $. The
statement of the theorem thus provides a complete account of the $K$-types
in $\mathbb{C}\left[ \overline{\mathcal{O}_{i}}\right] $, even when $i=n$
and $\Sigma =A_{n}$ or $a_{n}.$  \qed

\subsection{Remarks} 

\subsubsection{Connection with spherical orbits for symmetric pairs}
While our construction always yields a spherical nilpotent 
$K_{\frak{C}}$-orbit in $\frak{p}$ in the sense of \cite{king}, our construction falls short of producing all such
orbits. For example, we cannot obtain any spherical orbit whose closure is not 
a normal variety; for in our construction the ring of regular functions on an orbit
coincides with the ring of regular functions on its closure (well, except in the Hermitian
symmetric case where they nevertheless coincide up to a character of the center of $K$). To indicate which spherical 
nilpotent orbits are obtainable by our construction, we display in Appendix A how
the orbits we have constructed are situated within the Hasse diagrams depicting the closure 
relations (\cite{djokovic}) among the entire family of spherical nilpotent orbits 
as classifed by King (\cite{king}). 

From these diagrams one can understand the ``range'' of our construction 
in the following way: if $\Gamma=\left\{  \gamma_{1},\ldots,\gamma_{n}\right\}  $ is a
strongly orthogonal sequence of noncompact weights, the corresponding orbit closures
$
\overline{\mathcal{O}_{i}}$ 
are totally ordered by inclusion%
\[
\left\{  0\right\}  =\overline{\mathcal{O}_{0}}\subset \overline{\mathcal{O}_{1}}\subset
\overline{\mathcal{O}_{2}}\subset\cdots\subset\overline{\mathcal{O}_{n}}\quad.
\]
and are such that $\overline{\mathcal{O}_i} = \overline{\mathcal{O}_{i+1}} - \mathcal{O}_{i+1}$.
Such a sequence of orbits would correspond to a chain in the Hasse diagram in which there
are no branchings are encountered as one descends from the top of the chain to the trivial orbit. 
Nevertheless, in several cases our construction exhausts or nearly exhausts the 
set of spherical nilpotent orbits: for example, we get all the spherical orbits of
$SL ( n , \mathbb{R})$, $SL(n, \mathbb{H} )$, and $SO(3,p)$; and all but two 
orbits for $SO(2,p)$ and $Sp(p,q)$.

\subsubsection{Connection with unipotent representations}
We can now describe in a little more detail how we hope to attach unipotent
(in the sense of \cite{vogan}, Conjecture 12.1)
representations to these orbits. In \cite{sahi}, Siddhartha Sahi shows the
existence of a certain family of small unitary irreducible representations
of the conformal groups of simple non-Euclidean Jordan algebras.
Deliberately putting aside the Jordan theoretical underpinnings of these
representations, one can say that the essential representation-theoretical 
ingredients of Sahi's construction and analysis are:

\begin{itemize}
\item[(i)]  the circumstance that each representation is realized as a constituent 
of a non-unitary, degenerate, spherical principal series representation
\[ 
I(s) = Ind_{P=MAN}^{G}\left( 1\otimes e^{s\nu }\otimes 1\right)
\] 
whose 
associated variety is the closure of a single, multiplicity-free 
$K_{\mathbb{C}}$-orbit in $\frak{p}$; and

\item[(ii)]  the fact that there exists a $w\in N_{K}\left( \frak{a}\right) $ such that
both $wPw^{-1}=\overline{P}$ and $ad^{\ast }\left( w\right) \nu =-\nu $.
\end{itemize}

We remark that the second condition ensures that the underlying $\left( \frak{g}, K \right)$-module of $I(s)$
can be endowed with a $\frak{g}$-invariant
(but possibly indefinite or degenerate) Hermitian form (\cite{KZ}), while the first condition allows Sahi to
carry out an explicit
analysis of the action of $\frak{p}$ on $K$-types from which both signature
characters and reducibility conditions can be derived.

Using the results of \S 2 and \S 3, we can formulate a similar setup for any
connected semisimple Lie group, subsuming the situation of \cite{sahi}
in a uniform manner. Let $G$ be such a group, $\Gamma =\left[
\gamma _{1,}\ldots ,\gamma _{n}\right] $ a maximal sequence of strongly
orthogonal noncompact weights. As in Remark 2.1.2 we construct for each $%
i=1,\ldots ,n$, a normal $S$-triple $\left\{ x_{i},h_{i},y_{i}\right\} $
such that $x_{i}\in \frak{p}_{\gamma_{i}}$, $y_{i}=-\overline{x_{i}}\in \frak{p}%
_{-\gamma _{i}},h_{i}\in \frak{k}$. We set
\[
X_{n} =x_{1}+\cdots +x_{n} \quad , \quad
H_{n} =h_{1}+\cdots +h_{n} \quad , \quad
Y_{n} =y_{1}+\cdots +y_{n}
\]
as in (\ref{XHY}) and apply each of the (commuting) Cayley transforms (\ref{cayley}) 
successively to $\left\{X_{n},H_{n},Y_{n}\right\}$ to obtain a standard triple
\[
\left\{ X^{\prime}_{n}, H^{\prime}_{n}, Y^{\prime}_{n} \right\}
= 
\left\{
\frac{1}{2}\left( X_{n}+Y_{n}-iH_{n}\right) ,
-i\left( X_{n}-Y_{n}\right) ,
\frac{1}{2} \left( X_{n}+Y_{n}+iH_{n}\right) 
\right\}
\]
in real Lie algebra $\frak{g}_{\mathbb{R}}$ of $G$ such that the semisimple 
element $H^{\prime}_{n}$ is in $\frak{p}_{\mathbb{R}}$ and 
$\theta{X^{\prime}_{n}} = - Y^{\prime}_{n}$. Let
\begin{eqnarray*}
\frak{n} &=&\text{direct sum of positive eigenspaces of }ad\left(
H_{n}^{\prime }\right) \text{ in }\frak{g}_{\mathbb{R}} \quad , \\
\frak{l} &=&0\text{-eigenspace of }ad\left( H_{n}^{\prime }\right) \text{ in 
}\frak{g}_{\mathbb{R}} \quad , \\
\frak{a} &=&\left(\text{center of }\frak{l}\right) \cap \frak{p}_{\mathbb{R}} \quad , \\
\frak{m} &=&\text{orthogonal complement of }\frak{a}\text{ in }\frak{l}
\end{eqnarray*}
and set 
\begin{eqnarray*}
M &=&Z_{K}\left( \frak{a}\right) \exp \left( \frak{m}\right) \quad ,  \\
A &=&\exp \left( \frak{a}\right)  \quad , \\
N &=&\exp \left( \frak{n}\right) \quad . 
\end{eqnarray*}
Then $P=MAN$ is a (Langlands decomposition of a) parabolic subgroup of $G$.
Moreover, it happens that 
\begin{equation*}
span_{\mathbb{R}}\left( h_{1}^{\prime },\ldots ,h_{n}^{\prime
}\right) \subseteq \frak{a} \quad .
\end{equation*}
Now let $\nu $ be the element of the real dual space $\frak{a}^{\ast }$ of 
$\frak{a}$ such that 
\begin{equation*}
\nu \left( H\right) =B_{0}\left( H_{n}^{\prime },H\right) \qquad \forall \
H\in \frak{a} \quad ,
\end{equation*}
where $B_{0}\left( \cdot ,\cdot \right) $ is the Killing form on $\frak{g}_{\mathbb{R}}$
restricted to $\frak{a}$. Finally, we set 
\begin{equation*}
w=\exp \left( \frac{\pi }{2}\left( X_{n}^{\prime }-Y_{n}^{\prime }\right)
\right) \in K \quad .
\end{equation*}
It then happens that 
\[
w \in N_{K}\left( \frak{a}\right) \quad , \quad
wPw^{-1} =\overline{P} \quad , \quad
Ad^{\ast }\left( w\right)\nu  =-\nu \; .
\]
Thus, we have a natural means of attaching to our families of multiplicity-free
$K_{\mathbb{C}}$-orbits families of possibly unitarizable, degenerate principal 
series representations.

However, there is one more Jordan-theoretic device at play in
Sahi's paper; and that is a generalized Capelli operator $D_{1}$ that,
for a particular value of the parameter $s \in \mathbb{R}$, intertwines $I(s) $ 
with its Hermitian dual $I(-s)$.
And yet, even this Capelli operator should be characterizable in purely
representation-theoretic terms as per \cite{Boe}.  In fact, we conjecture here there 
is there is a quasi-invariant differential operator on $C^{\infty}\left( \overline{\frak{n}} \right)$
corresponding to the Cayley transform of a lowest weight vector of the irreducible representation
homogeneous summand of $S^{n}\left( \frak{p}\right) $ of highest weight $%
\gamma _{1}+\cdots +\gamma _{n}$ that intertwines a $I(s)$ with its Hermitian dual.  
In a subsequent paper, we hope to confirm
(or correct) this conjecture and to extend the methods of \cite{KS} and \cite{sahi}
to study of the
reducibility and signature characters of the subrepresentations $I_{P}\left(
s \right) $ associated with the orbits $K\cdot Y_{i}$.

\section{Dimension and Degree of $\overline{\mathcal{O}_{i}}$}

Let $R=\bigoplus_{n=0}^{\infty} R_{n}$ be the polynomial ring $\mathbb{C}\left[  x_{1},\ldots,x_{n}\right]  $
regarded as a graded commutative ring (the grading by degree of homogenity). If $I\subset R$ is
homogeneous ideal then $M=R/I$ is a graded $R$-module:%
\[
 M=\bigoplus_{n=0}^{\infty} M_{n} \quad \text{where} \quad
M_{n}\equiv R_{n}/\left(  R_{n}\cap I\right) \; .
\]
Let $Y$ be the corresponding affine variety. 
By a theorem of Hilbert and Serre, there is a unique
polynomial $p_{Y}\left(  t\right)  $ such that%
\[
p_{Y}\left(  t\right)  = \sum_{k=0}^{t}\dim_{\mathbb{C}}M_{k}%
 \qquad\text{for all
sufficiently large }t \quad .%
\]
$p_{Y}\left(  t\right)  $ is called the Hilbert polynomial of $Y$. When one
writes%
\[
p_{Y}\left(  t\right)  =ct^{d}+ \left( \text{terms of order } t^{d-1}\right)
\]
the degree of the leading term is the \emph{dimension} (often by definition)
of $Y$ and the number%
\[
D=\frac{c}{d!}%
\]
is the \emph{degree} of $Y$. It turns out that $D$ is always an
integer and corresponds to the number of points where the projectivization of
$Y$ meets a generic $\left(  n-d-1\right)  $-dimensional linear subspace of
$\mathbb{P}^{n-1}$.

Now consider a $K_{\mathbb{C}}$-orbit  $\mathcal{O}_{i}$ associated to a
sequence  $\Gamma_{i} = \left\{  \gamma_{1},\ldots,\gamma_{i}\right\}  $ of strongly
orthogonal noncompact weights and  the ring of regular functions $\mathbb{C}%
\left[\overline{\mathcal{O}_{i}}\right] $ on its closure. We shall assume for ease of exposition that $i<n=|\Gamma |$.
The cases when $i=n$ can be handled similarly; but not so uniformly.
It follows from the proof of Theorem 3.11
that a given $K$-type $\lambda=a_{1}\gamma_{1}+\cdots+a_{i}\gamma_{i}$ in
$\mathbb{C}\left[  \mathcal{O}_{i}\right]  $ is generated by the action of
$\frak{k}$ on
\[
\phi_{\lambda}=\left(  \phi_{\omega_{1}}\right)  ^{a_{1}-a_{2}}\left(
\phi_{\omega_{2}}\right)  ^{a_{2}-a_{3}}\cdots\left(  \phi_{\omega_{i}%
}\right)  ^{a_{i}}%
\]
which is a homogeneous polynomial of degree (taking $a_{j}=0$ for
$j>i$)%
\begin{equation*}
\sum_{j=1}^{i}j\left(  a_{j}-a_{j+1}\right) 
 =\sum_{j=1}^{i}a_{j} \quad .%
\end{equation*}
It follows that%
\[
p_{\overline{\mathcal{O}_{i}}}\left(  t\right)  =\sum_{\ell=0}^{t}\dim
M_{\ell}=\sum_{\mathbf{\lambda}\in\Lambda_{t}}\dim V_{\lambda
}%
\]
where
\[
\Lambda_{t}=\left\{  a_{1}\gamma_{1}+\cdots+a_{i}\gamma_{i}\mid a_{1}%
,\ldots,a_{i}\in\mathbb{N}\;,\; a_{1}\geq a_{2}\geq\cdots\geq a_{i}%
\geq0\;,\;\sum_{j=1}^{i}a_{i}\leq t\right\} \;.
\]
Applying the Weyl dimension formula,  we obtain%
\[
p_{\overline{\mathcal{O}_{i}}}\left(  t\right)  =\sum_{\lambda\in\Lambda_{t}%
}\left(  \prod_{\alpha\in\Delta^{+}\left(  \frak{t};\frak{k}\right)  }%
\frac{\left\langle \lambda+\rho_{K},\alpha\right\rangle }{\left\langle
\rho_{K},\alpha\right\rangle }\right)  \quad.
\]
Now note that the factors $\left\langle
\lambda+\rho_{K},\alpha\right\rangle /\left\langle \alpha,\rho_{k}%
\right\rangle $ either reduce to factors of $1$ (when $\alpha\perp\lambda
$)\ or contribute factors of the form $\left\langle \lambda,\alpha
\right\rangle /\left\langle \alpha,\rho_{K}\right\rangle $ to the leading term
of $p_{\overline{O_{i}}}$.   
We thus need only account for the roots that have components along
$\gamma_{1},\ldots,\gamma_{i}$.  Let
\begin{equation*}
\Delta_{i}^{+}   =\left\{  \alpha\in\Delta^{+}\left(  \frak{t};\frak{k}%
\right)  \mid\left\langle \alpha,\gamma_{j}\right\rangle \neq0\text{ for some
}j\in\left\{  1,\ldots,i\right\}  \right\}
\end{equation*}
The leading term of the Hilbert polynomial for $p_{\overline{\mathcal{O}}_{i}}\left(  t\right)$ 
is thus
\begin{equation*}
LT\left(  p_{\overline{\mathcal{O}}_{i}}\right) 
 =\left(  \prod_{\alpha\in\Delta_{i}^{+}}\frac{1}{\left\langle \rho
_{K},\alpha\right\rangle }\right)  LT\left(  \sum_{\lambda\in\Lambda_{t}%
}\left(  \prod_{\alpha\in\Delta_{i}^{+}}\left\langle \lambda,\alpha
\right\rangle \right)  \right)
\end{equation*}

To compute the products
$
\prod_{\alpha\in\Delta_{i}^{+}}\left\langle \lambda,\alpha
\right\rangle  $, 
we just need to know $\Sigma$, the restricted root systems with
multiplicities associated with $\left\{ \gamma_{1}, \ldots , \gamma_{n} \right\}$.
Since each restricted root will be of one of the types $a_{n}$,
$A_{n},\ldots,d_{n}$,and since the roots of each type share a common multiplicity, 
we can write 
\begin{equation*}
\prod_{\alpha \in \Delta _{i}^{+}}\left\langle \lambda ,\alpha \right\rangle
=\left( \prod_{\alpha \in \left( a_{n}\right) _{i}^{+}}\left\langle \lambda
,\alpha \right\rangle \right) ^{m_{a}}\left( \prod_{\alpha \in \left(
A_{n}\right) _{i}^{+}}\left\langle \lambda ,\alpha \right\rangle \right)
^{m_{A}}\cdots \cdots \left( \prod_{\alpha \in \left( d_{n}\right)
_{i}^{+}}\left\langle \lambda ,\alpha \right\rangle \right) ^{m_{d}}
\end{equation*}
where 
\begin{align*}
\left( a_{n}\right) _{i}^{+}& =\left\{ \frac{1}{2}\gamma _{j}-\frac{1}{2}%
\gamma _{k}\mid 1\leq j<k\leq i\right\} \cup \left\{ \frac{1}{2}\gamma _{j}-%
\frac{1}{2}\gamma _{k}\mid 1\leq j\leq i<k\leq n\right\} \quad , \\
\left( A_{n}\right) _{i}^{+}& =\left\{ \gamma _{j}-\gamma _{k}\mid 1\leq
j<k\leq i\right\} \cup \left\{ \gamma _{j}-\gamma _{k}\mid 1\leq j\leq
i<k\leq n\right\} \quad , \\
\left( b_{n}\right) _{i}^{+}& =\left\{ \frac{1}{2}\gamma _{j}\mid 1\leq
j\leq i\right\} \quad , \\
\left( C_{n}\right) _{i}^{+}& =\left\{ \gamma _{j}\mid 1\leq j\leq i\right\}
\quad , \\
\left( a_{n}\right) _{i}^{+}& =\left\{ \frac{1}{2}\gamma _{j}\pm \frac{1}{2}%
\gamma _{k}\mid 1\leq j<k\leq i\right\} \cup \left\{ \frac{1}{2}\gamma
_{j}\pm \frac{1}{2}\gamma _{k}\mid 1\leq j\leq i<k\leq n\right\} \quad .
\end{align*}
We thus have 
\begin{eqnarray*}
\left( \prod_{\alpha \in \left( a_{n}\right) _{i}^{+}}\left\langle \lambda
,\alpha \right\rangle \right) ^{m_{a}} &=&\left( \frac{1}{2}\right)
^{m_{a}i\left( 2n-i-1\right) /2}\left( \prod_{1\leq j<k\leq i}\left(
a_{j}-a_{k}\right) ^{m_{a}}\right) \left( \prod_{1\leq j\leq i}\left(
a_{j}\right) ^{m_{a}\left( n-i\right) }\right) \quad , \\
\left( \prod_{\alpha \in \left( A_{n}\right) _{i}^{+}}\left\langle \lambda
,\alpha \right\rangle \right) ^{m_{A}} &=&\left( \prod_{1\leq j\leq i}\left(
a_{j}\right) ^{m_{A}\left( n-i\right) }\right) \left( \prod_{1\leq j<k\leq
i}\left( a_{j}-a_{k}\right) ^{m_{A}}\right) \quad , \\
\left( \prod_{\alpha \in \left( b_{n}\right) _{i}^{+}}\left\langle \lambda
,\alpha \right\rangle \right) ^{m_{b}} &=&\left( \frac{1}{2}\right)
^{im_{b}}\left( \prod_{1\leq j\leq i}\left( a_{j}\right) ^{m_{b}}\right)
\quad , \\
\left( \prod_{\alpha \in \left( C_{n}\right) _{i}^{+}}\left\langle \lambda
,\alpha \right\rangle \right) ^{m_{C}} &=&\left( \prod_{1\leq j\leq i}\left(
a_{j}\right) ^{m_{C}}\right) \quad , \\
\left( \prod_{\alpha \in \left( d_{n}\right) _{i}^{+}}\left\langle \lambda
,\alpha \right\rangle \right) ^{m_{d}} &=&\left( \frac{1}{2}\right)
^{m_{d}i\left( 2n-i-1\right) }\left( \prod_{1\leq j<k\leq i}\left(
a_{j}+a_{k}\right) ^{m_{d}}\right)  \\
&&\quad \times \left( \prod_{1\leq j<k\leq i}\left( a_{j}-a_{k}\right)
^{m_{d}}\right) \left( \prod_{1\leq j\leq i}\left( a_{j}\right) ^{2\left(
n-i\right) m_{d}}\right) \quad .
\end{eqnarray*}
And so we get
\begin{align*}
\left( \prod_{\alpha \in \Delta _{i}^{+}}\left\langle \lambda ,\alpha
\right\rangle \right) & =\left( \frac{1}{2}\right) ^{\frac{1}{2}i\left(
2n-i-1\right) \left( m_{a}+2m_{d}\right) +im_{b}}\left( \prod_{1\leq j\leq
i}a_{j}\right) ^{m_{b}+m_{C}+\left( n-i\right) \left(
m_{a}+m_{A}+2m_{d}\right) } \\
& \times \left( \prod_{1\leq j<k\leq i}\left( a_{j}-a_{k}\right) \right)
^{m_{a}+m_{A}}\left( \prod_{1\leq j<k\leq i}\left(
a_{j}^{2}-a_{k}^{2}\right) \right) ^{m_{d}}\quad .
\end{align*}
In summary, 

\begin{Lemma}
\begin{equation*}
LT\left( p_{\overline{\mathcal{O}}_{i}}\right) =c_{i}\ LT\left( \sum_{\lambda
\in \Lambda _{t}}\left( \prod_{1\leq j\leq i}a_{j}\right) ^{q}\left(
\prod_{1\leq j<k\leq i}\left( a_{j}-a_{k}\right) \right) ^{r}\left(
\prod_{1\leq j<k\leq i}\left( a_{j}^{2}-a_{k}^{2}\right) \right) ^{s}\right) 
\end{equation*}
where
\begin{align*}
c_{i}& =\frac{1}{\prod_{\alpha \in \Delta _{i}^{+}}\left\langle \rho _{K},\alpha
\right\rangle }\left( \frac{1}{2}\right) ^{\frac{1}{2}i\left( 2n-i-1\right)
\left( m_{a}+2m_{d}\right) +im_{b}}\quad , \\
q& =\left( n-i\right) \left( m_{a}+m_{A}+2m_{d}\right) +\left(
m_{b}+m_{C}\right) \quad , \\
r& =m_{a}+m_{A}\quad , \\
s& =m_{d}\quad .
\end{align*}
\end{Lemma}

Next we observe that the sum over $\Lambda_{t}$ can be carried out as an
iterated sum of the form
\[
\sum_{\lambda\in\Lambda_{t}}\left(  \cdots\right)  =\sum_{a_{1}=0}^{t}%
\sum_{a_{2}=0}^{\min\left(  a_{1},t-a_{1}\right)  }\cdots\sum_{a_{i}=0}%
^{\min\left(  a_{i-1},t-a_{1}-\cdots-a_{i-i}\right)  }\left(  \cdots\right)
\]
and that the quantity to be summed,%
\[
F\left(  \mathbf{a}\right)  \equiv\left(  \prod_{j=1}^{i}\left(  a_{j}\right)
^{q}\right)  \left(  \prod_{1\leq j<k\leq i}\left(  a_{j}%
-a_{k}\right)  ^{r}\left(  a_{j}^{2}-a_{k}^{2}\right)  ^{s}\right)
\]
is homogeneous of degree%
\[
deg(F) = qi + i(i-1)(r+2s)/2 
\]
in the variables $a_{1},\ldots,a_{i}$.

\begin{Lemma}
Suppose $\Omega_{t}\subset\mathbb{N}^{n}$ is a region of the form%
\[
\Omega_{t}=\left\{  \mathbf{a}\in\mathbb{N}^{n}\mid0\leq a_{1}\leq
t\;,\;0\leq a_{2}\leq\phi_{2}\left(  t,a_{1}\right)  \;,\;
\cdots\;,0\leq a_{n}\leq\phi_{n}\left(  t,a_{1},\ldots,a_{n-1}\right)
\right\}  \;,
\]
where each $\phi_{i}$ is a homogeneous linear function of its arguments. Then
for large $t$%
\begin{eqnarray*}
\lefteqn{\sum_{\mathbf{a}\in\Omega_{t}}a_{1}^{m_{1}}\cdots a_{n}^{m_{n}} =}  \\ 
&&\int_{0}%
^{t}\int_{0}^{\phi_{2}\left(  t,x_{1}\right)  }\cdots\int_{0}^{\phi_{n}\left(
t,x_{1},\ldots,x_{n-1}\right)  }x_{1}^{m_{1}}\cdots x_{n}^{m_{n}}dx_{n}\cdots
dx_{1}+\text{ lower order terms}\quad.
\end{eqnarray*}
\end{Lemma}

\emph{Proof.} From Faulhaber's formula \cite{faulhaber}
\[
\sum_{i=0}^{t}t^{p}=\frac{1}{p+1}\sum_{k=1}^{p+1}\left(  -1\right)
^{\delta_{k,p}}\left(
\begin{array}
[c]{c}%
p+1\\
k
\end{array}
\right)  B_{p+1-k}t^{k}%
\]
(where $\delta_{k,p}$ is the Kronecker delta symbol, $\left(
\begin{array}
[c]{c}%
a\\
b
\end{array}
\right)  $ is the usual binomial coefficients, and $B_{q}$ is the $q^{th}$
Bernoulli number) one sees that
\begin{align*}
S_{p}\left(  t\right)   &  \equiv\sum_{i=0}^{t}i^{p}=\frac{1}{p+1}%
t^{p+1}+\frac{1}{2}t^{p}-\frac{p}{12}t^{p-1}+\cdots\\
&  =\int_{0}^{t}x^{p}dx+\text{ lower order terms} \quad .%
\end{align*}
The result now follows from an easy computation and inductive argument. \qed

\begin{Remark}
Note that
\[
P\left(  t\right)  =\int_{0}^{t}\int_{0}^{\phi_{2}\left(  t,x_{1}\right)
}\cdots\int_{0}^{\phi_{n}\left(  t,x_{1},\ldots,x_{n-1}\right)  }x_{1}^{m_{1}%
}\cdots x_{n}^{m_{n}}dx_{n}\cdots dx_{1}%
\]
is a monomial of degree $m_{1}+\cdots+m_{n}+n$ in $t$. Its leading
coefficient is thus%
\[
P\left(  1\right)  =\int_{0}^{1}\int_{0}^{\phi_{2}\left(  1,x_{1}\right)
}\cdots\int_{0}^{\phi_{n}\left(  1,x_{1},\ldots,x_{n-1}\right)  }x_{1}^{m_{1}%
}\cdots x_{n}^{m_{n}}dx_{n}\cdots dx_{1} \quad .%
\]
Thus if we set%
\[
R_{t}=\left\{  \mathbf{x}\in\mathbb{R}^{n}\mid0\leq x_{1}\leq 1 , 0\leq
x_{2}\leq\phi_{2}\left(  1,x_{1}\right)  , \cdots ,0\leq x_{n}\leq\phi
_{n}\left(  1,x_{1},\ldots,x_{n-1}\right)  \right\}
\]
we have%
\[
\sum_{\mathbf{a}\in\Omega_{t}}a_{1}^{m_{1}}\cdots a_{n}^{m_{n}}\approx\left(
\int_{R_{1}}x_{1}^{m_{1}}\cdots x_{n}^{m_{n}}d^{n}x\right)  t^{m_{1}%
+\cdots+m_{n}+n}\quad.
\]
\end{Remark}

\begin{Proposition}
Let
\[
\mathcal{S}_{n,t}=\left\{  \mathbf{x}\in\mathbb{R}^{n}\mid x_{1}\geq x_{2}%
\geq\cdots\geq x_{n}\geq0\ ,\ \sum_{i=1}^{n}x_{i}\leq t\right\}
\]
and%
\[
\Lambda_{t}=\mathcal{S}_{n,t}\cap\mathbb{N}^{n}%
\]
If $F\left(  x_{1},\ldots,x_{n}\right)  $ is homogeneous of degree $d$, then
\[
\sum_{\mathbf{a}\in\Lambda_{t}}F\left(\mathbf{a}\right)  =\left(
\int_{\mathcal{S}_{n,1}}F\left( \mathbf{x}\right)  dx_{n}\cdots
dx_{1}\right)  t^{d+n}+\text{lower order terms} \quad .%
\]
\end{Proposition}

\emph{Proof.} We can decompose the sum over
$\Lambda_{t}$ into a sum of sums%
\[
\sum_{\Lambda_{t}}F\left(  a_{1},\cdots,a_{n}\right)  =\sum_{i=1}^{N}%
\sum_{\Lambda_{t,i}}F\left(  a_{1},\cdots,a_{n}\right)
\]
where each region $\Lambda_{t,i}$ is a region of the elementary form
considered in Lemma (4.1). One can then apply the lemma and the subsequent
remark to get%
\[
\sum_{\Lambda_{t}}F\left(  a_{1},\cdots,a_{n}\right)  \approx\sum_{i=1}%
^{N}\left(  \int_{R_{1},i}F\left(  x_{1},\cdots,x_{n}\right)  d^{n}x\right)
t^{d+n}+\text{ lower order terms}%
\]
and then reassemble the region $\mathcal{S}_{n,1}$ from the regions
$R_{1,i}$ to obtain the desired result. \qed

Applying the preceding proposition to our formula for the Hilbert polynomial
$p_{\overline{\mathcal{O}_{i}}}\left(  t\right)  $, we can conclude:

\begin{Theorem}
Let $G$ be a semisimple Lie group, $\left\{  \gamma_{1},\ldots,\gamma
_{n}\right\}  $ a maximal sequence of strongly orthogonal noncompact weights,
$\Sigma$ the corresponding restricted root system (specified as in (\ref{0})), and
let $\mathcal{O}_{i}$ be the multiplicity-free $K_{\mathcal{C}}$-orbit
associated to a subsequence $\left\{  \gamma_{1},\ldots,\gamma_{i}\right\}  $.
The dimension of $\overline{\mathcal{O}_{i}}$ is given by%
\[
\dim\left(  \mathcal{O}_{i}\right)  = i(q+1) + i(i-1)(r+2s)/2 
\]
and its degree is given by%
\begin{align*}
Deg\left( \overline{\mathcal{O}_{i}}\right)   &  =\frac{c_{i}}{\dim\left(  \mathcal{O}%
_{i}\right)  !}\left(  \prod_{\alpha\in\Delta^{+}\left(  \frak{t}%
;\frak{k}\right)  }\frac{1}{\left\langle \rho_{K},\alpha\right\rangle
}\right) \\
&  \times\int_{\mathcal{S}_{i}}\left(  \prod_{j=1}^{i}x_{i}\right)
^{q}\left(  \prod_{1\leq j<k\leq i}\left(  x_{j}-x_{k}%
\right)  \right)  ^{r}\left(  \prod_{1\leq j<k\leq i}\left(
x_{j}^{2}-x_{k}^{2}\right)  \right)  ^{s}d^{i}x
\end{align*}
where $\mathcal{S}_{i}$ is the domain%
\[
\mathcal{S}_{i}=\left\{  x\in\mathbb{R}^{i}\mid x_{1}\geq x_{2}\geq\cdots\geq
x_{i}\geq0\quad,\quad\sum_{j=1}^{i}x_{j}\leq1\right\} \quad 
\]
and the constants $c_i$, $q$, $r$ and $s$ are as in Lemma 4.1 and Table 2 below.
\end{Theorem}

\begin{table}\caption{}
\begin{equation*}
\begin{array}{|l|c|c|c|c|}%
\hline
 G  &  \dim\mathcal{O}_{i} & q & r & s \\ \hline
SL\left(  n,\mathbb{R}\right)   & i\left(  2\left[  \frac{n}{2}\right]
-i\right)   & 2\left(  \left[  \frac{n}{2}\right]  -i\right)   & 0 & 1 \\ \hline 
 SL\left(  n,\mathbb{H}\right)   & 4i\left(  2\left[  \frac{n}{2}\right]
-i\right)   & 8\left(  \left[  \frac{n}{2}\right]  -i\right)  +3 & 0 & 4 \\ \hline 
 SU\left(  p,q\right)   & i\left(  p+q-i\right)   & p+q-2i & 2 & 0 \\ \hline 
 SO\left(  2,q\right)  \ ,\ q>2 & i\left(  i\left(  2-q\right)  +q-4\right)
/2 & \left(  q-2\right)  \left(  2-i\right)   & q-2 & 0 \\ \hline 
 SO\left(  p,q\right)  \ I\ ,\ 2<p\leq q & i\left(  p+q-2i-1\right) &
p+q-4i & 0 & 2\\ \hline 
 SO^{\ast}(2n) & i\left(  8\left[  \frac{n}{2}\right]  -4i-1\right) &
8\left(  \left[  \frac{n}{2}\right]  -i\right)  +2 & 0 & 4 \\ \hline 
 Sp\left(  n,\mathbb{R}\right)   & i\left(  2n-i+1\right)  /2 & n-i & 1 & 0 \\ \hline 
 Sp\left(  p,q\right)   & 2i\left(  p+q-i+1\right)   & 2\left(
p+q-2i+1\right)   & 0 & 2 \\
\hline 
\end{array}
\end{equation*}
\end{table}

\begin{Remark}
From the data tabulated above, one finds that the formula for the degree of
$\overline{\mathcal{O}_{i}}$ is either of the form%
\begin{equation}
\int_{S_{i}}\left(  \prod_{j=1}x_{j}\right)  ^{q}\left(  \prod_{1\leq j<k\leq
i}\left(  x_{j}-x_{k}\right)  \right)  ^{r}d^{i}x\label{4}%
\end{equation}
(which happens only when $G/K$ is Hermitian symmetric), or%
\begin{equation}
\int_{S_{i}}\left(  \prod_{j=1}x_{j}\right)  ^{q}\left(  \prod_{1\leq j<k\leq
i}\left(  x_{j}^{2}-x_{k}^{2}\right)  \right)  ^{s}d^{i}x\label{5} \quad .%
\end{equation}
The first form is very much akin to the famous Selberg integral \cite{selberg}%
\[
S_{n,r,s,t}=\frac{1}{n!}\int_{\left[  0,1\right]  ^{n}}\left(  \prod_{i=1}%
^{n}x_{i}\right)  ^{s}\left(  \prod_{i=1}^{n}\left(  1-x_{i}\right)  \right)
^{t}\left(  \prod_{1\leq i<j\leq n}\left|  x_{i}-x_{j}\right|  \right)
^{r}d^{n}x \quad . 
\]
Indeed, by setting $t=0$ and making a change of variables (\cite{macdonald}, pg. 286 )
(\ref{4}) can be explicitly evaluated using Selberg's formula. The more
generic case (\ref{5}), however, seems to be lacking an explicit evaluation.
We do note, however, that Nishiyama, Ochiai and Zhu \cite{noz} encountered 
and evaluated certain integrals of the form (\ref{5}) in their study of theta liftings of 
nilpotent orbits.  In \cite{binegar} we provide several other methods of 
evaluating integrals of the general form (\ref{5}). 

\end{Remark}

\appendix

\section{Closure Relations for Spherical Nilpotent Orbits of Classical Real Linear Groups}

To indicate exactly which spherical orbits are constructible by our sequences of
strongly orthogonal noncompact weights, we display below
the closure relations (\cite{ohta} , \cite{djokovic}) for the spherical orbits of classical real linear
Lie groups(\cite{king}); or rather those cases for which we've identified a nice pattern (the closure
diagrams of $SU\left( p,q \right)$ and $SO \left( p,q \right)$ get rather complicated
as $p$ and $q$ increase). The double lines in the diagram indicate the simple chains of
spherical nilpotent orbits closures corresponding to sequences of strongly orthogonal noncompact 
weights (cf. Remark 2.1.3.). In the Hermitian symmetric cases we indicate both the chains lying 
in $\frak{p}_{+}$ and those lying in
$\frak{p}_{-}$.   Our notation for the orbits is somewhere between that 
of \cite{king} and \cite{djokovic}.
Briefly, as in \cite{king} we indicate particular orbits by expressions of the form 
$(\pm n_{1})^{m_{1}} (\pm n_{2})^{m_{2}}  \cdots (\pm n_{k})^{m_{k}}$, where a factor of the form
$(\pm n_{i})^{m_{i}}$ indicates the occurance of a signed
a row of alternating '+' and '-' signs, of length $n_{i}$, beginning with a $\pm$ sign,
and occuring with multiplicity $m_{i}$. However, Djokovic's algorithm makes use of unsigned
rows (actually, unsigned ``genes'') rather than rows that are more commonly represented 
as even signed rows; we indicate such an unsigned row of length $n$ occuring with multiplicity $m$
by a factor of the form $\left(n\right)^{m}$. Thus, for example,

\def\hunit{pt}
\def\vunit{pt}
\def\place#1 #2 #3 {\rlap{\kern#1\hunit
        \raise#2\vunit\hbox{#3}}}

\newdimen\hoogte    \hoogte=12pt    
\newdimen\breedte   \breedte=14pt   
\newdimen\dikte     \dikte=0.5pt    
\def\beginYoung{
       \begingroup
       \def\vr{\vrule height0.8\hoogte width\dikte depth 0.2\hoogte}
       \def\fbox##1{\vbox{\offinterlineskip
                    \hrule height\dikte
                    \hbox to \breedte{\vr\hfill##1\hfill\vr}
                    \hrule height\dikte}}
       \vbox\bgroup \offinterlineskip \tabskip=-\dikte \lineskip=-\dikte
            \halign\bgroup &\fbox{##\unskip}\unskip  \crcr }
\def\endYoung{\egroup\egroup\endgroup}

\begin{center}
\hbox{
\place 70 40 {$(+3)^{2}(2)(+1)^{2}(-1) \quad \sim \quad$}
\place 190 0 {\beginYoung
     +  &  $-$ &  +\cr     
     +  &  $-$ & +   \cr
     ~  &  ~     \cr
     +           \cr
     + \cr
     $-$ \cr
\endYoung}
\place 260 40 {.}
} 
\end{center}

\begin{figure}[h]\caption{$SL\left( n,\mathbb{R} \right)$}
\small
\[
\xymatrix@C=3pt@R=10pt{
\mathcal{O}_{(2)^{\lbrack n/2 \rbrack}}^{I}\ar@{=>}[dr]
&&\mathcal{O}_{(2)^{\lbrack n/2 \rbrack}}^{II}\ar@{=>}[dl]\\
&\qquad\mathcal{O}_{(2)^{\lbrack n/2 \rbrack -1}(1)^{2}}\ar@{=>}[d] \\
&\qquad\mathcal{O}_{(2)^{\lbrack n/2 \rbrack -2}(1)^{4}}\ar@{=>}[d] \\
&\vdots \ar@{=>}[d]\\
&\mathcal{O}_{(2)(1)^{n-2}}\ar@{=>}[d] \\
&\mathcal{O}_{(1)^{n}} \\
&n \text{ even}
}
\hspace{.4in}
\xymatrix@C=3pt@R=10pt{
\quad \mathcal{O}_{(2)^{\lbrack n/2 \rbrack}}\ar@{=>}[d]\\
\qquad\mathcal{O}_{(2)^{\lbrack n/2 \rbrack -1}(1)^{2}}\ar@{=>}[d] \\
\qquad\mathcal{O}_{(2)^{\lbrack n/2 \rbrack -2}(1)^{4}}\ar@{=>}[d] \\
\vdots \ar@{=>}[d]\\
\mathcal{O}_{(2)(1)^{n-2}}\ar@{=>}[d] \\
\mathcal{O}_{(1)^{n}}\\
n \text{ odd} 
}
\]
\normalsize
\end{figure}

\begin{figure}[h]\caption{$SU\left( 2, q \right)$}
\small
\[
\xymatrix@C=3pt@R=10pt{
&&\mathcal{O}_{\left(-3\right)^{2}\left( -1 \right)^{q-4}}\ar[drr]\ar[dll]\\
\mathcal{O}_{\left( -3\right)  \left(  -2\right)\left( -1 \right)^{q-3}  }\ar [dd] \ar [drr]&&&&
\mathcal{O}_{\left(  -3\right)  \left( +2\right) \left( -1 \right)^{q-3}  } \ar [dll] \ar [dd] \\
&
&\mathcal{O}_{\left( -3\right)  \left( +1\right)  \left(  -1\right)^{q-1}  } \ar[d]
&\mathcal{O}_{\left(  +3\right)  \left(  -1\right) ^{q-1}} \ar[dl]  \\
\mathcal{O}_{\left(  -2\right) ^{2}\left(  -1\right)^{q-2}  }\ar\ar@{=>}[dr]&&
\mathcal{O}_{\left(  +2\right)  \left(  -2\right) \left(  -1\right)^{q-2}  }\ar[dl]\ar[dr]
&&\mathcal{O}_{\left(  +2\right)  ^{2} \left(  -1\right)^{q-2}  }\ar\ar@{=>}[dl]\\
&\mathcal{O}_{\left(  -2 \right)  \left(  +1\right) \left( -1\right)^{q-1}} \ar\ar@{=>}[dr]
&&\mathcal{O}_{\left(  +2\right)  \left(  +1\right)  \left(  -1\right)^{q-1}  }\ar@{=>}[dl]\\
&&\mathcal{O}_{\left(  +1\right)  ^{2}\left(  -1\right)  ^{q}}
}
\]
\normalsize
\end{figure}

\begin{figure}[h]\caption{$SL\left( n,\mathbb{H} \right)$}
\small
\[
\xymatrix@C=3pt@R=10pt{
\quad \mathcal{O}_{(2)^{\lbrack n/2 \rbrack}}\ar@{=>}[d]\\
\qquad\mathcal{O}_{(2)^{\lbrack n/2 \rbrack -1}(1)^{2}}\ar@{=>}[d] \\
\vdots \ar@{=>}[d]\\
\mathcal{O}_{(2)(1)^{n-2}}\ar@{=>}[d] \\
\mathcal{O}_{(1)^{n}}\\
n \text{ even} 
}
\hspace{.8in}
\xymatrix@C=3pt@R=10pt{
\quad \mathcal{O}_{(2)^{\lbrack n/2 \rbrack}(1)}\ar@{=>}[d]\\
\qquad\mathcal{O}_{(2)^{\lbrack n/2 \rbrack -1}(1)^{3}}\ar@{=>}[d] \\
\vdots \ar@{=>}[d]\\
\mathcal{O}_{(2)(1)^{n-2}}\ar@{=>}[d] \\
\mathcal{O}_{(1)^{n}}\\
n \text{ odd} 
}
\]
\normalsize
\end{figure}

\begin{figure}[h]\caption{$SO\left( 2, p \right) \quad ; \quad p>4$}
\small
\[
\xymatrix@C=3pt@R=10pt{
&&\mathcal{O}_{(-3)^{2}(-1)^{p-4}} \ar[d] \\
\mathcal{O}_{(+3)(-1)^{p-1}}^{I}\ar@{=>}[dr]
&&\mathcal{O}_{(-3)(+1)(-1)^{p-2}} \ar[dl] \ar[dr] 
&&\mathcal{O}_{(+3)(-1)^{p-1}}^{II}\ar@{=>}[dl]\\
&\mathcal{O}_{(+2)(-2)(-1)^{p-2}}^{I} \ar@{=>}[dr]
&&\mathcal{O}_{(+2)(-2)(-1)^{p-2}}^{II} \ar@{=>}[dl] \\
&&\mathcal{O}_{(+1)^{2}(-1)^{p}} 
}
\]
\normalsize
\end{figure}

\begin{figure}[h]\caption{$SO^{\ast} \left( 2n \right)$}
\small
\[
\xymatrix@C=1pt@R=10pt{
\mathcal{O}_{\lbrack \frac{n}{2} \rbrack,0}\ar@{=>}[dr]
&&\mathcal{O}_{\lbrack \frac{n}{2} \rbrack -1,1}\ar[dl]\ar[dr]
&\cdots
&\cdots
&\cdots
&\cdots
&\cdots
&\mathcal{O}_{1,\lbrack \frac{n}{2} \rbrack -1}\ar[dr]\ar[dl]
&&\mathcal{O}_{0,\lbrack \frac{n}{2} \rbrack} \ar@{=>}[dl]  \\
&\mathcal{O}_{\lbrack \frac{n}{2} \rbrack -1,0} \ar@{:>}[ddrr]
&\cdots
&~
&&
&
&~
&
&\mathcal{O}_{0,\lbrack \frac{n}{2} \rbrack -1} \ar@{:>}[ddll] \\
&&&&&\mathcal{O}_{(3)(1)^{n-3}}\ar[d]&&&\\
&&&\mathcal{O}_{2,0}\qquad\ar@{=>}[dr]
&&\mathcal{O}_{1,1}\ar[dl]\ar[dr]
&&\qquad\mathcal{O}_{0,2}\ar@{=>}[dl]
\\
&&&&\mathcal{O}_{1,0}\ar@{=>}[dr]
&&\mathcal{O}_{0,1}\ar@{=>}[dl]\\
&&&&&\mathcal{O}_{0,0}
}
\]
\normalsize
where
\small
$\mathcal{O}_{r,s} = \mathcal{O}_{(+2)^{r}(-2)^{s}(1)^{n-2r-2s}}$
\normalsize
\end{figure}

\begin{figure}[h]\caption{$Sp\left( n, \mathbb{R} \right)$}
\small
\[
\xymatrix@C=3pt@R=10pt{
\mathcal{O}_{n,0}\ar@{=>}[dr]
&&\mathcal{O}_{n-1,1}\ar[dl]\ar@{.>}[ddrr]
&\cdots
&\cdots
&\cdots
&\mathcal{O}_{1,n-1}\ar[dr]\ar@{.>}[ddll]
&&\mathcal{O}_{0,n} \ar@{=>}[dl]  \\
&\mathcal{O}_{n-1,0} \ar@{:>}[ddrr]
&\cdots
&
&\cdots
&
&
&\mathcal{O}_{0,n-1} \ar@{:>}[ddll] \\
&&&&\mathcal{O}_{1,1}\ar[dl]\ar[dr]
\\
&&&\mathcal{O}_{1,0}\ar@{=>}[dr]
&&\mathcal{O}_{0,1}\ar@{=>}[dl]\\
&&&&\mathcal{O}_{0,0}
}
\]
\normalsize
where
$\mathcal{O}_{r,s} = \mathcal{O}_{(+2)^{r}(-2)^{s}(+1)^{n-r-s}(-1)^{n-r-s}}$
\normalsize
\end{figure}

\begin{figure}[h]\caption{$Sp(p,q) \quad p \le q$}
\small
\[
\xymatrix@C=3pt@R=10pt{
&\mathcal{O}_{(+2)^{p}}\ar@{=>}[d]\\
&\mathcal{O}_{(+2)^{p-1}(+1)(-1)}\ar@{=>}@{:>}[dd]\\
\\
\mathcal{O}_{(+3)(+2)(+1)^{p-3}(-1)^{q-2}}\ar[d]\ar[dr]
&\mathcal{O}_{(+2)^{3}(+1)^{p-3}(-1)^{q-3}}\ar@{=>}[d]
&\mathcal{O}_{(-3)(+2)(+1)^{p-2}(-1)^{q-3}}\ar[dl]\ar[d]
\\
\mathcal{O}_{(+3)(+1)^{p-2}(-1)^{q-1}}\ar@{=>}[dr]
&\mathcal{O}_{(+2)^{2}(+1)^{p-2}(-1)^{q-2}}\ar@{=>}[d]
&\mathcal{O}_{(-3)(+1)^{p-1}(-1)^{q-2}}\ar@{=>}[dl]
\\
&\mathcal{O}_{(+2)(+1)^{p-1}(-1)^{q-1}}\ar@{=>}[d]
\\
&\mathcal{O}_{(+1)^{p}(-1)^{q}}
}
\]
\normalsize
\end{figure}

\clearpage
\normalsize


\begin{thebibliography}{99}
\bibitem[B]{binegar}B. Binegar, \emph{On the evaluation of some Selberg-like integrals}, to appear.
(\texttt{arXiv:math.RT/0608301}, 
http:$\backslash\backslash$arxiv.org$\backslash$abs$\backslash$math.RT$\backslash$0608301). 

\bibitem[Bo]{Boe} B. Boe, \emph{Homomorphisms between generalized Verma modules}, 
Trans. Am. Math. Soc., \textbf{288} (1985), 791-799.

\bibitem[Bour] {bourbaki}N. Bourbaki, Groupes et Alg\`{e}bres de Lie, Chapitres
4,5, et 6, Masson, Paris (1981).

\bibitem[Br1]  {Br1}M. Brion, \emph{Groupe de Picard et nobres caracteristiques des
varietes spheriques}, Duke Math. J. \textbf{58} (1989), no. 2, 397-424.

\bibitem[Br2] {Br2}M. Brion, \emph{Spherical Varieties: an introduction}, in
Topological methods in algebraic transformation groups (New Brunswick, NJ
1988),  Prog. Math. \textbf{80} (1989), Birkhauser Boston, Boston, MA, 11-26.

\bibitem[D] {djokovic}D. Djokovic, \emph{Closures of Conjugacy Classes in Classical Real 
Linear Lie Groups}, in Algebra, Carbondale, 1980, 
Lecture Notes in Math. \textbf{848}, Springer-Verlag, New York,
1980, 63-88.

\bibitem[F] {faulhaber}J. Faulhaber, Academiae Algebrae (1631).

\bibitem[Ka] {kac}V. Kac, \emph{Some Remarks on Nilpotent Orbits}, J. of
Algebra \textbf{64} (1980), 190-213.

\bibitem[KO] {kato-ochiai}S. Kato and H. Ochiai, \emph{The degrees of orbits of
the multiplicity free actions}, in ``Nilpotent Orbits, Associated Cycles and
Whittaker Models for Highest Weight Representations'', Ast\'{e}risque
\textbf{273}, Soci\`{e}t\`{e} Math. France (2001).

\bibitem[Ki] {king}D. King, \emph{Classification of Spherical Nilpotent Orbits
in Complex Symmetric Space}, J. Lie Theory \textbf{14} (2004), 339-370.

\bibitem[Ko] {kostant}B. Kostant, \emph{Lie Group Representations on
Polynomial Rings}, Am. J. Math. \textbf{86} (1963), 327-402.

\bibitem[KR] {kostant-rallis} B. Kostant and S. Rallis, \emph{Orbits and
representations associated with symmetric spaces}, Amer. J. Math. \textbf{93}
(1971), 753-809.

\bibitem[KS]{KS}  B. Kostant and S. Sahi, \emph{Jordan Algebras and
Capelli Identities}, Inventiones Math., \textbf{112} (1993), 657-664.

\bibitem[KV] {kimelfeld-vinberg}B.N. Kimel'fel'd and E.B. Vinberg,
\emph{Homogeneous Domains on Flag Manifolds and Spherical Subgroups of
Semisimple Lie Groups}, Functional Anal. Appl., \textbf{12} (1978), 12-19.

\bibitem[Ku] {kumar}S. Kumar, \emph{Proof of the Parthasarathy-Ranga
Rao-Varadarajan conjecture}, Inv. Math. \textbf{93} (1988), 117-130.

\bibitem[KY] {KY}H. Kaji and O. Yasukura, \emph{Secant varieties of adjoint
varieties:\ orbit decomposition}, J. Algebra \textbf{227} (2000), no. 1, 26-44.

\bibitem[KZ] {KZ} A. Knapp and G. Zuckerman, \emph{Classification theorems
for representations of semisimple Lie groups}, Non-Commutative Harmonic Analysis,
Springer Verlag Lec. Notes in Math. \textbf{587} (1977), 138-159.

\bibitem[M] {macdonald} I.G. Macdonald, Symmetric Functions and Hall
Polynomials, Clarendon\ Press, Oxford, 1995.

\bibitem[MRS] {mrs} I. Muller, H. Rubenthaler, and G. Schiffmann,
\emph{Structure des espaces pr\'ehomogen\`enes associ\'e \`a certaines 
alg\`ebres de Lie gradu\'ees}, Math. Ann. \textbf{274} (1986), no. 1, 95-123.

\bibitem[N1]  {N1}K. Nishiyama, \emph{Multiplicity-free actions and the geometry of
nilpotent orbits}, Math. Ann. \textbf{318} (2000), no. 4, \ 777-793.

\bibitem[N] {N2}K. Nishiyama, \emph{Classification of spherical nilpotent
orbits for $U(p,q)$}, J. Math. Kyoto Univ. \textbf{44} (2004), 203-215.

\bibitem[NO] {NO} K. Nishiyama and H. Ochiai, \emph{Bernstein degree of
singular unitary highest weight representations of the metaplectic group},
Proc. Japan. Acad., \textbf{75} (1999), 9-11. 

\bibitem[NOT] {NOT} K. Nishiyami, H. Ochiai, and K. Taniguchi, 
\emph{Bernstein degree and associated cycles of Harish-Chandra modules - 
Hermitian symmetric case}, in ``Nilpotent Orbits, Associated Cycles and
Whittaker Models for Highest Weight Representations'', Ast\'{e}risque
\textbf{273}, Soci\`{e}t\`{e} Math. France (2001).

\bibitem[NOZ] {noz} K. Nishiyama, H. Ochiai, and C. Zhu, \emph{Theta liftings
of nilpotent orbits for symmetric pairs},
Trans. Amer. Math. Soc.  358  (2006),  no. 6, 2713--2734.

\bibitem[O] {ohta} T. Ohta, \emph{The closures of nilpotent orbits in classical
symmetric pairs and their singularities}, Tohoku Math. J. (2) \textbf{43}
(1991), no. 2, 161-211.

\bibitem[Pa] {panyushev} D. Panyushev, \emph{On spherical nilpotent orbits and beyond}, 
Ann. Inst. Fourier (Grenoble) \textbf{49} (1999), no.5, 1453-1476.

\bibitem[Sa1] {sahi0} S. Sahi, \emph{Unitary Representations on the Shilov Boundary of a 
Symmetric Tube Domain}, Contemp. Math. \textbf{145} (1993), 275-286.

\bibitem[Sa2] {sahi} S. Sahi, \emph{Jordan algebras and degenerate principal
series}, J. reine angew. Math. \textbf{462} (1995), 1-18.


\bibitem[Se] {selberg} A. Selberg, \emph{Bemerkninger om et multipelt
integral}, Norsk. Mat. Tidsskr, \textbf{26} (1944), 71-78.

\bibitem[Sev] {servedio} F.J. Servedio, \emph{Prehomogeneous vector spaces and
varieties}, Trans. Am. Math. \textbf{176} (1973), 421-444.

\bibitem[Vo] {vogan} D. Vogan, \emph{Associated Varieties and Unipotent
Representations}, in ``Harmonic Analysis on Reductive Groups'', Prog. in
Math., \textbf{101}, Birkh\"auser Boston, Boston, MA, 1991, 315-388.
\end{thebibliography}
\end{document}